\RequirePackage{fix-cm}
\documentclass[smallextended]{svjour3}       % onecolumn (second format)
\smartqed  % flush right qed marks, e.g. at end of proof
\usepackage{graphicx}
\usepackage{etex,bm}
\usepackage{amsmath,amssymb,mathrsfs,mathtools}
\usepackage{amsfonts}
\usepackage[margin=1in]{geometry}
\usepackage{multirow,booktabs}
\usepackage{empheq}
\usepackage{array}
\usepackage{subcaption}
\usepackage[colorlinks,citecolor=red,urlcolor=blue,bookmarks=false,hypertexnames=true]{hyperref}
\begin{document}
\title{Feasibility Conditions of Robust Portfolio Solutions under Single and Combined Uncertainties%\thanks{Grants or other notes
%about the article that should go on the front page should be
%placed here. General acknowledgments should be placed at the end of the article.}
}
%\subtitle{Do you have a subtitle?\\ If so, write it here}

%\titlerunning{Short form of title}        % if too long for running head

\author{Pulak Swain         \and
        Akshay Kumar Ojha %etc.
}

%\authorrunning{Short form of author list} % if too long for running head

\author{Pulak Swain* \thanks{*Corresponding author} \and Akshay Kumar Ojha}
\institute{Pulak Swain \at
	School of Basic Sciences, Indian Institute of Technology Bhubaneswar \\
	%Fax: +123-45-678910\\
	\email{ps28@iitbbs.ac.in}           %  \\
	%             \emph{Present address:} of F. Author  %  if needed
	\and
	Akshay Kumar Ojha \at
	School of Basic Sciences, Indian Institute of Technology Bhubaneswar\\
	\email{akojha@iitbbs.ac.in}
}

\date{Received: date / Accepted: date}
% The correct dates will be entered by the editor

\maketitle
\begin{abstract}
In this paper, we derive the feasibility conditions for the robust counterparts of the uncertain Markowitz model. Our study is based on ellipsoidal, box, polyhedral uncertainty sets and also the uncertainty sets obtained from their combinations. We write the uncertain Markowitz problem in a quadratically constrained convex quadratic programming form and the uncertain sets are converted into their quadratic forms for the derivation. The feasibility conditions of robust solutions which we obtained are in the form of positive semi-definite matrices. The study can be useful for deriving the robust counterparts of the combined uncertainties, which in general is difficult to find out.\\
% \PACS{PACS code1 \and PACS code2 \and more}
% \subclass{MSC code1 \and MSC code2 \and more}
\end{abstract}

\section{Introduction} \label{sec4.1}
\par Markowitz's theory \cite{Markowitz J 1952} for portfolio allocation has been an evolution in the world of finance. The theory is based on the process of allocating total wealth among different assets in the best possible way with the help of the optimization model. Risk and reward are the two major components in such problems, which are conflicting in nature. An investor always wants to invest in a manner that he would get maximum reward from the portfolio by taking the minimum risk. In Markowitz's model the expected return and varaince of return are respectively the reward and risks associated with the portfolio.  Later on, some downside risk measures \cite{Konno 2002,Estrada 2007} such as semivariance, semiabsolute deviation, below target risk, value at risk and conditional value at risk has been used, mainly when the asset return distributions are asymmetrical. Moreover, some more models such as capital asset pricing model, sharpe ratio model were also used to study portfolio optimization \cite{Sharpe 1964,Dowd 2000}. However the Markowitz model has the highest popularity among the researchers due to its simplicity and transparency while solving the problem. The most important aspect of this model is the impact on portfolio diversification by the number of assets and their covariance relationships \cite{Markowitz 1959}. With the passage of time, many studies emerge in the field of portfolio optimization. Meghwani and Thakur \cite{Meghwani 2017} introduced some practical constraints such as, cardinality constraint, bound constraint, pre-assignment constraint, which can be added in the problem as per the needs of an investor.  In recent years many researchers make the use of heuristic methods to
study portfolio optimization. Cura \cite{Cura 2009} applied particle swarm optimization approach to a portfolio optimization problem and compared the results with those of the genetic algorithm, simulated annealing and tabu search approaches. Branke et. al. \cite{Branke 2009} studied the portfolio optimization with an envelope-based multi-objective evolutionary algorithm. Meghwani and Thakur \cite{Meghwani 2018} proposed a novel heuristic based repair algorithm for multi-objective portfolio optimization with practical constraints.
\par One of the biggest challenges in portfolio optimization is that the portfolio models are highly sensitive to the uncertainty factor and the model may give a bad decision if the uncertainty is not handled well. The main reason for uncertainty in these models is the use of historical data for getting future decisions. The term "uncertainty" refers to the ambiguity in constant coefficients in the optimization problem. This is due to the lack of knowledge in exact value at the time when the problem is being solved. Previously the uncertain optimization problems were being dealt with help of dynamic programming \cite{Bellman 1954},
fuzzy optimization \cite{Lodwick 2010}, sensitivty analysis \cite{Christopher 2002}. However, in the last two decades robust optimization approach has been widely used by the researchers due
to its ability to produce completely ”immunized against uncertainty” solutions. In this approach the problem of worst case realization of uncertain parameters is solved so that the solution becomes valid for any realization of uncertainty. One of the primary assumptions in robust optimization is that the uncertain parameters vary in a prespecified set $\mathscr{U}$.  The robust optimization approach aims to choose the best solution among those “immunized” against data uncertainty, that is, the solutions should remain feasible for all the realization of data within $\mathscr{U}$. The "worst-case-oriented" philosophy makes
it natural to quantify the quality of a robust feasible solution by the guaranteed value of the original objective. ElGhaoui and Lebret \cite{El Ghaoui 1997} first studied the robust solutions for the uncertain least-squares problems, and Ben-Tal \& Nemirovski \cite{Ben-Tal 1998} laid the foundation on robust convex optimization. El-Ghaoui et al. \cite{El Ghaoui 1998} studied uncertain semidefinite
problems. They provided sufficient conditions which guarantee that the robust solution is unique and continuous for the unperturbed problem's data. Ben-Tal and Nemirovski \cite{Ben-Tal 1999} studied the robust counterparts for several types of convex optimization problems for ellipsoidal uncertainty.
The study showed that when the uncertainty sets for a linear constraint are ellipsoids, the robust formulation turns out to be a conic quadratic problem. Now a days robust optimization has a wide range
of applications \cite{Gabrel 2014} in various fields such as inventory \& logistics, finance, revenue management, queueing networks \& game theory, machine learning, energy systems etc.
\par After the introduction of robust optimization, many researchers in the financial research field also make use of this in their problems. Goldfarb and Iyengar \cite{Goldfarb 2003} used robust optimization to a factor model for representing return. T\"{u}t\"{u}nc\"{u} and Koenig \cite{Tutuncu 2004} formulated the robust counterparts for bound-based uncertainty sets and also they illustrated the robust efficient frontier. Natarajan et. al. \cite{Natarajan 2009} used robust optimization to generalize different portfolio risk measures and showed that coherent risk measures can be constructed from the uncertainty sets in robust optimization models. Fabozzi et al. \cite{Fabozzi 2010} applied the robust optimization approach for mean-variance, value at risk, and conditional value at risk models. Fliege and
Werner \cite{Fliege 2014} applied the robust multiobjective optimization in the portfolio mean-variance problem. Kim et al. \cite{Kim 2018} did an analysis on the performance of robust portfolio for U.S. equity portfolios and their study highlighted the effectiveness of robust optimization for controlling uncertainty in efficient investments. Sehgal and Mehra \cite{Sehgal 2019} proposed robust reward-risk ratio portfolio models using Omega, semi-mean absolute deviation ratio, and weighted stable tail adjusted return ratio. Their results show that the robust reward-risk models outperformed their conventional counterpart models. Ismail and Pham \cite{Ismail 2019} studied a robust continuous-time Markowitz model with an uncertain covariance matrix of risky assets. They also obtained a lower bound of Sharpe ratio of any robust efficient portfolio strategy.
\par Since the robust optimization framework uses the worst-case parameter approach, the feasibility of a robust solution is very much important. Because more often the solutions are infeasible for the worst case realizations of uncertain parameters. In the existing literature, there has not been any proper study regarding the feasibility conditions for the robust portfolio problems and this is the prime objective of this paper. The main contribution of the paper is as follows. Robust feasibility conditions of Markowitz portfolio model are derived for several uncertainties such as box, ellipsoidal, polyhedral, etc. Further, the study is extended for the combined uncertainties, that is the intersections between any of the two aforementioned uncertain sets. The uncertain constraints as well as the uncertain sets are first expressed in quadratic form and then the feasibility conditions are obtained in the form of positive semidefinite matrices. Since the robust counterparts of combined uncertainties have not been explored much, the present study can be used to derive the robust counterparts of uncertain Markowitz model under these combined uncertainties.
\par The paper is organized as follows: The introduction part is followed by some preliminary ideas on the Markowitz portfolio model, robust optimization and the uncertain sets in Section \ref{sec4.2}. Section \ref{sec4.3} presents the feasibility conditions for the robust solutions of the Markowitz model under single as well as combined uncertainty sets. Finally, some concluding remarks have been provided in Section \ref{sec4.4}.
\section{Preliminaries} \label{sec4.2}
\subsection{Markowitz Portfolio Model}
Suppose there are $n$ number of assets $S_1, S_2, \dots, S_n$ with their respective returns at time $t$ be given by $r_{1t}, r_{2t}, \dots, r_{nt}$. Markowitz's portfolio model is based on two input parameters-- expected portfolio return (acts as the reward) and the variance of portfolio return (acts as the risk factor). To calculate these, we first need to find the expected return of each asset $(\mu_i)$ and the covariance of return between each pair of assets $(\sigma_{ij})$, which are given by:
\begin{align*}
\begin{aligned}	
&\mu_i=E(r_i)=\frac{1}{T}\sum\limits_{t=1}^{T}r_{it}, \qquad \sigma_{ij}=E[(r_i-\mu_i)(r_j-\mu_j)]=\frac{1}{T}\sum\limits_{t=1}^{T}(r_{it}-\mu_i)(r_{jt}-\mu_j)\\
&\textrm{for } i=1, 2,\dots,n \textrm{ and } j=1, 2, \dots,n.
\end{aligned}
\end{align*}
Let $\bm{\mu}$ be the $(n \times 1)$ vector of all $\mu_i$'s and $\bm{\Sigma}$ be the $(n \times n)$ matrix of all $\sigma_{ij}$'s. The aim is to form a portfolio which will give our desired return with a minimum risk associated with it. Let the weight given to asset $S_i$ be $x_i$. Then the expected return and variance of return for the portfolio are respectively given by,
\begin{align*}
{\mu_P} = \sum_i \mu_i x_i= \bm{\mu}^\top\bm{x}, \quad 	{\sigma_P}^2 = \sum_{i,j} {\sigma_{ij}x_ix_j} = \bm{x}^\top\bm{\Sigma} \bm{x}
\end{align*}
where $\bm x$ be the $(n \times 1)$ decision vector containing the weights $x_i$'s.\\
Markowitz's portfolio model minimizes the variance of portfolio return at a given level of expected portfolio return (say $\tau$). The mathematical formulation of the model is given by,
\begin{align} 
	\begin{aligned}	
		&\min_{\bm{x}} && \frac{1}{2} \bm{x}^\top\bm{\Sigma} \bm{x} \\
		&\text{s.t.:} && \bm{\mu}^\top \bm{x} \geq \tau, \quad	\bm{e}^\top\bm{x}=1, \quad \bm{x}\geq \bm{0}
	\end{aligned}
\label{eq4.2}
\end{align}
 and $\bm e$ is the $(n \times 1)$ vector containing $1$'s, that is, $\bm e=(1 \  1 \ \dots \ 1)^\top$.
%Equivalently if we represent an ${(n \times 1)}$ vector $\bm{\mu}$ as the expected return vector containing $\mu_i$'s and  an $(n \times n)$ matrix $\bm{\Sigma}$ as the covariance matrix containing each $\sigma_{ij}$, then the model (1) can  be interpreted in vector form as:
%\begin{align} 
%	\begin{aligned}	
%		&\min_{\bm{x}} && \frac{1}{2} \bm{x}^\top\bm{\Sigma} \bm{x} \\
%		&\text{s.t.:} && \bm{\mu}^\top \bm{x} \geq \tau, \qquad	\bm{1}^\top\bm{x}=1, \qquad \bm{x}\geq \bm{0}
%	\end{aligned}
%\end{align}
%where $\bm x$ be the $(n \times 1)$ decision vector containing the weights $x_i$'s and $\bm 1$ is the $(n \times 1)$ vector containing $1$'s.
\subsection{Uncertainty in Optimization Problems and the Robust Optimization Approach}
An \textit{uncertain optimization problem} is defined as the optimization problem in which the constant coefficients are not exactly known; the only thing we know is that they are close to some nominal values. The general form of an uncertain optimization problem is given by,
\begin{align}
	\begin{aligned}	
		&\min_{\bm x} && f(\bm{x},\bm{u}) \\
		&\textrm{s.t.:} && c(\bm{x},\bm{u}) \leq 0, \quad \bm u  \in \mathscr{U} 
	\end{aligned}
	\label{eq4.2a}
\end{align}
where $\bm x$ is the vector of decision variables, $\bm u$ is the vector of uncertain parameters and $\mathscr{U}$ is the uncertain set.\\
The uncertain set can be affinely parameterized by a perturbation vector in a given perturbation set $\mathcal{Z}$ as follows:
\begin{align}
\begin{aligned}	
\mathscr{U}=\left\{\bm{u} : \bm{u}= \bm{u}^{(0)}+\sum_{j=1}^{n}\zeta_j \bm{u}^{(j)}, \ \bm{\zeta} = (\zeta_1, \zeta_2, \dots, \zeta_n)\ \in \mathcal{Z} \subset \mathbb{R}^n \right\}
\end{aligned}
\label{eq4.4}
\end{align}
where $\bm{u}^{(0)}$ is the nominal vector of the uncertain parameters, $\bm{u}^{(j)}$ are the basic shifts, and $\bm{\zeta}$ is the perturbation vector.\\
For an example, suppose we have three uncertain parameters $u_1$, $u_2$ and $u_3$ in our problem. Let the nominal values are respectively $u_1^{(0)}=3$,  $u_2^{(0)}=4$ and $u_3^{(0)}=5$. And let the values of $u_1$, $u_2$ and $u_3$ can perturb upto $0.02$, $0.04$ and $0.03$ respectively. So we can write the uncertainty set as,
\begin{align*}
\begin{aligned}	
\mathscr{U}=\left\{ \begin{bmatrix} 
u_1\\u_2\\u_3
\end{bmatrix} : \begin{bmatrix} 
3\\4\\5
\end{bmatrix}+ \zeta_1 \begin{bmatrix} 
0.02\\0\\0 \end{bmatrix} + \zeta_2 \begin{bmatrix} 
0\\0.04\\0 \end{bmatrix} + \zeta_3 \begin{bmatrix} 
0\\0\\0.03 \end{bmatrix},\ \zeta_1, \zeta_2, \zeta_3 \in [-1, 1] \right\}
\end{aligned}
\end{align*}
Some important concepts and definitions of robust optimization relevant to the present study are enlisted as follows:
\begin{definition} [Robust Feasible Solution] \label{def4.1} \cite{Ben-Tal 2009}
	A vector $\bm x$ is said to be the robust feasible solution of the uncertain problem \eqref{eq4.2a} if it satisfies the uncertain constraints for all realizations of the uncertain set $\mathscr{U}$, that is, if $\bm x$ satisfies
	\begin{align*}
	\begin{aligned}	
	c(\bm{x},\bm{u}) \leq 0, \quad \forall \bm u  \in \mathscr{U} 
	\end{aligned}
	\end{align*}
\end{definition}
\begin{definition} [Robust Value] \label{def4.2} \cite{Ben-Tal 2009}
	For a given candidate solution $\bm{x}$, the robust value $\widehat{f}(\bm{x})$ of the objective in problem \eqref{eq4.2a} is the largest value of $f(\bm{x},\bm{u})$ over all realizations of the data from the uncertain set, that is
	\begin{align*}
	\begin{aligned}	
	\widehat{f}(\bm{x})= \max_{\bm u  \in \mathscr{U}} f(\bm{x},\bm{u})
	\end{aligned}
	\end{align*}
\end{definition}
It is to be noted that if the problem \eqref{eq4.2a} was a maximization problem, then the robust value would have been $\widehat{f}(\bm{x})= \min_{\bm u  \in \mathscr{U}} f(\bm{x},\bm{u})$.
\begin{definition} [Robust Counterpart] \label{def4.3} \cite{Ben-Tal 2009}
	The robust counterpart of the uncertain problem \eqref{eq4.2a} is the optimization problem 
	\begin{align}	
	\begin{aligned} \label{eq4.2b}
	& \min_{\bm x} && \left \{ \max_{\bm u \in \mathscr{U}}	f(\bm x, \bm u) \right \}\\
	&\text{s.t.:} && c(\bm x, \bm u) \leq \bm 0, \quad \forall \bm u  \in \mathscr{U}
	\end{aligned}
	\end{align}
	of minimizing the robust value of the objective over all robust feasible solutions to the uncertain problem.
\end{definition}
In other words, the robust counterpart is the optimization problem with the worst case realization of uncertain parameters.
\begin{definition} [Robust Optimal Solution] \label{def4.4} \cite{Ben-Tal 2009}
	The optimal solution of the robust counterpart problem \eqref{eq4.2b} or is said to be the robust optimal solution of the uncertain problem \eqref{eq4.2a}.
\end{definition}
However, sometimes the robust counterpart problems are very difficult to solve. So the tractability of robust counterparts should be taken into consideration when solving uncertain problems.
\subsection{Several Uncertainties} \label{sec4.2.3}
Robust optimization assumes the uncertain sets in the form of bounded convex sets. Some of the most used uncertainties are given as,\\
(i) Box Uncertainty: \begin{align*}\mathscr{U}^{(B)}=\left \{\bm{\zeta}:||\bm{\zeta}||_{\infty} \leq \delta_B \right \} = \left \{{\zeta_j}:\sup_j |{\zeta_j}| \leq \delta_B \right \} \end{align*}
(ii) Ellipsoidal Uncertainty: \begin{align*}\mathscr{U}^{(E)}=\left \{\bm{\zeta}:||\bm{\zeta}||_{2} \leq \delta_E \right \}= \left \{{\zeta_j}:\sum_j({\zeta_j})^2 \leq (\delta_E)^2 \right \}\end{align*}
(iii) Polyhedral Uncertainty:\begin{align*}\mathscr{U}^{(P)}=\left \{\bm{\zeta}:||\bm{\zeta}||_{1} \leq \delta_P \right \}=\left \{{\zeta_j}:\sum_j|{\zeta_j}| \leq \delta_P \right \}\end{align*}
where $\delta_B$, $\delta_E$ and $\delta_P$ are respectively the radii of box, eliipsoidal and polyhedral uncertain sets.\\
Li and Floudas \cite{Li 2012} used the combinations of the above uncertainties. Since the intersection of two bounded convex sets is also bounded and convex, these combinations can be used as the uncertain sets. The combined sets are enlisted as follows.\\
(i) "Box $\cap$ Ellipsoidal" Uncertainty Set:\begin{align*}\mathscr{U}^{(BE)}=\left \{\bm{\zeta}:||\bm{\zeta}||_{\infty} \leq \delta_B,\ ||\bm{\zeta}||_{2} \leq \delta_E \right \}=\left \{{\zeta_j}:\sup_j |{\zeta_j}| \leq \delta_B,\ \sum_j({\zeta_j})^2 \leq (\delta_E)^2 \right \}\end{align*}
(ii) "Box $\cap$ Polyhedral" Uncertainty Set:\begin{align*}\mathscr{U}^{(BP)}=\left \{\bm{\zeta}:||\bm{\zeta}||_{\infty} \leq \delta_B,\ ||\bm{\zeta}||_{1} \leq \delta_P \right \}=\left \{{\zeta_j}:\sup_j |{\zeta_j}| \leq \delta_B,\ \sum_j|{\zeta_j}| \leq \delta_P \right \}\end{align*}
(iii) "Ellipsoidal $\cap$ Polyhedral" Uncertainty Set:\begin{align*}\mathscr{U}^{(EP)}=\left \{\bm{\zeta}:||\bm{\zeta}||_{2} \leq \delta_E,\ ||\bm{\zeta}||_{1} \leq \delta_P \right \}=\left \{{\zeta_j}:\sum_j({\zeta_j})^2 \leq (\delta_E)^2,\ \sum_j|{\zeta_j}| \leq \delta_P \right \}\end{align*}
%(iv) "Box $\cap$ Ellipsoidal $\cap$ Polyhedral" Uncertainty Set:\begin{align*}\mathscr{U}^{(BEP)}=\left \{\bm{\zeta}:||\bm{\zeta}||_{\infty} \leq \delta_B,\ ||\bm{\zeta}||_{2} \leq \delta_E,\ ||\bm{\zeta}||_{1} \leq \delta_P \right \}=\left \{{\zeta_i}:\max_i |{\zeta_i}| \leq \delta_B,\ \sum_i({\zeta_i})^2 \leq (\delta_E)^2,\ \sum_i|{\zeta_i}| \leq \delta_P \right \}\end{align*}
\section{Feasibility Conditions of Robust Markowitz Model with Single and Combined Uncertainties} \label{sec4.3}
 The study of Ben-Tal and Nemirovski \cite{Ben-Tal 1999} on convex optimization problems includes an uncertain linear program with quadratic constraints given by,
\begin{align}
\begin{aligned} \label{eq4.3}
&\min_x && \bm c ^\top \bm x  \\
&\textrm{s.t.:} && - \bm x^\top [\bm A^i]^\top \bm A^i \bm x+ 2 [\bm b^i]^\top \bm x+\gamma^i > 0, \quad \bm (\bm A^i,\bm b^i,\gamma^i) \in \mathscr{U}_i\\
& && i=1, 2, \dots, m
\end{aligned}
\end{align} 
where $\mathscr{U}_i$ is the $i^{th}$ component of a bounded ellipsoid given by,
\begin{align}
\begin{aligned} \label{eq4.3aa}	
&\mathscr{U}_i=\left \{(\bm A^i,\bm b^i,\gamma^i)=(\bm A^{i(0)},\bm b^{i(0)},\gamma^{i(0)})+\sum_{j=1}^{k} \zeta_j (\bm A^{i(j)},\bm b^{i(j)},\gamma^{i(j)}) : \zeta_1^2+\zeta_2^2+\dots+\zeta_k^2 \leq 1 \right \}
\end{aligned}
\end{align}
The robust counterpart of the problem obtained from the study is equivalent to an explicit semidefinite programming problem (SDP). The result is given as follows.
\begin{theorem} \cite{Ben-Tal 1999} \label{thm4.0}
The robust counterpart of the uncertain problem \eqref{eq4.3} associated with the uncertain set $\mathscr{U}=\mathscr{U}_1 \times \mathscr{U}_2 \times \dots \times \mathscr{U}_m$, $\mathscr{U}_i$'s being given by \eqref{eq4.3aa}, is equivalent to the SDP problem
\begin{align}
\begin{aligned}	\label{eq4.3a}
&\min_x && \bm c ^\top \bm x  \\
&\textrm{w.r.t.} && \bm x \in \mathbb{R}^n, \ \lambda^1, \dots, \lambda^n \in \mathbb{R}\\
& \textrm{s.t.:}&& \begin{small}\left [\begin{array} {c|c c c c|c}
\gamma^{i(0)}+2 \bm x^\top\bm b^{i(0)}-\lambda^i & \dfrac{\gamma^{i(1)}}{2}+\bm x^\top \bm b^{i(1)}  & \dfrac{\gamma^{i(1)}}{2}+\bm x^\top \bm b^{i(1)} & \dots&\dfrac{\gamma^{i(1)}}{2}+\bm x^\top \bm b^{i(1)} & [\bm A^{i(1)}\bm x]^\top\\
\hline
\dfrac{\gamma^{i(1)}}{2}+\bm x^\top \bm b^{i(1)} & \lambda^i & & & & [\bm A^{i(1)}\bm x]^\top\\
\dfrac{\gamma^{i(2)}}{2}+\bm x^\top \bm b^{i(2)} & & \lambda^i & & & [\bm A^{i(2)}\bm x]^\top\\
\vdots & & &\ddots& & \vdots \\
\dfrac{\gamma^{i(k)}}{2}+\bm x^\top \bm b^{i(k)} & & & & \lambda^i & [\bm A^{i(k)}\bm x]^\top\\
\hline
\bm A^{i(0)}\bm x & \bm A^{i(1)}\bm x & \bm A^{i(2)}\bm x & \dots & \bm A^{i(k)}\bm x & I_{l_i}
\end{array}\right ]\end{small} \succeq 0\\
& && \text{for } i=1,2,\dots,m \text{ and $I_{l_i}$ is the identity matrix of order $l_i$}.
\end{aligned}
\end{align} 
\end{theorem}
During the derivation of the robust counterpart, the ellipsoidal uncertain sets $\mathscr{U}_i$'s, and the uncertain constraints are derived in quadratic form. Consequently the robust counterpart \eqref{eq4.3a} is obtained with the help of the well known S-lemma. 
\begin{lemma}[S-Lemma] \cite{Ben-Tal 2009} \label{lemma4.1} Let $\bm A$, $\bm B$ be symmetric matrices of the same size such that, there exists $\bm{\bar z}$ satisfying $\bm{\bar z}^\top \bm A \bm{\bar z} \geq 0$. Then the implication
	\begin{align*}
	\begin{aligned}	
	\bm z^\top \bm A \bm z \geq 0 \implies \bm z^\top \bm B \bm z \geq 0
	\end{aligned}
	\end{align*}
	holds true if and only if
	\begin{align*}
	\begin{aligned}	
	\exists \  \lambda \geq 0 \quad : \quad \bm B \succeq \lambda \bm A
	\end{aligned}
	\end{align*}
\end{lemma}
This approach can be used for deriving the robust feasibility conditions of uncertain quadratic constraint(s) for ellipsoidal uncertainty. In fact, the feasibility conditions for the other uncertainties can be obtained, if they are expressed in quadratic forms. We aim to derive the feasibility conditions of uncertain  Markowitz portfolio model under several single and combined uncertainties given in Section \ref{sec4.2.3}.
\subsection{Construction of the Uncertain Markowitz Portfolio Problem}
The uncertainty in the covariance of asset returns does not affect the optimal solution much as compared to the uncertainty in expected asset returns does \cite{Pulak 2021}. So we assume that the uncertainty is only due to the assets' expected returns.
Then the the uncertain Markowitz portfolio model is given by,
\begin{align}
\begin{aligned}	\label{eq4.5}
&\min_{\bm x} && \frac{1}{2} \bm{x}^\top\bm{\Sigma} \bm{x} \\
&\textrm{s.t.:} && \bm{\mu}^\top \bm{x} \geq \tau, \quad \bm{\mu} \in \mathscr{U}_{\bm{\mu}} \\
& &&	\bm{e}^\top\bm{x}=1, \quad \bm{x}\geq \bm{0}
\end{aligned}
\end{align}
where the uncertain set corresponding to $\mathscr{U}_{\mu}$ is given by,
\begin{align*}
\begin{aligned}	
\mathscr{U}_{\mu}=\left\{\bm{\mu} : \bm{\mu}= \bm{\mu}^{(0)}+\sum_{j=1}^{n}\zeta_j \bm{\mu}^{(j)}, \ \bm{\zeta} = (\zeta_1, \zeta_2, \dots, \zeta_n)^\top \in \mathcal{Z} \right\}
\end{aligned}
\end{align*}
The problem \eqref{eq4.4} becomes a QCQP problem by converting the uncertain constraint into a quadratic constraint, which is given by,
\begin{align}
\begin{aligned}	\label{eq4.7}
&\min_{\bm x} && \frac{1}{2} \bm{x}^\top\bm{\Sigma} \bm{x} \\
&\textrm{s.t.:} && \bm{x}^\top \bm{\mu} \bm{\mu}^\top \bm{x} - \tau^2 \geq 0, \quad \bm{\mu} \in \mathscr{U}_{\bm{\mu}} \\
& && \bm{e}^\top\bm{x}=1, \quad \bm{x}\geq \bm{0}.
\end{aligned}
\end{align}
Next we derive the robust feasibility conditions for single as well as combined uncertainties. 
\subsection{Feasibility Conditions for Single Uncertainties}
Since the robust counterpart result discussed earlier is for ellipsoidal uncertainty, we begin with the ellipsoidal uncertainty case.
\subsubsection{Ellipsoidal Uncertainty}
When the uncertain set in the problem \eqref{eq4.7} is a unit ellipsoid, it is represented by,
\begin{align}
\begin{aligned}	\label{eq4.6}
\mathscr{U}_{\mu}^{(E)}=\left\{\bm{\mu} : \bm{\mu}= \bm{\mu}^{(0)}+\sum_{j=1}^{n}\zeta_j \bm{\mu}^{(j)}, \ ||\bm{\zeta}||_2 \leq 1 \right\}
\end{aligned}
\end{align}
Then any solution $\bm{x}$ is robust feasible for the uncertain Markowitz potfolio  model with ellipsoidal uncertainty, if and only if we have,
\begin{align}
\begin{aligned}	\label{eq4.8}
& \bm{x}^\top \left[\bm{\mu}^{(0)}+\sum_{j=1}^{n}\zeta_j \bm{\mu}^{(j)}\right] \left[\bm{\mu}^{(0)}+\sum_{j=1}^{n}\zeta_j \bm{\mu}^{(j)}\right]^\top\bm{x} - \tau^2 \geq 0, \quad \forall \bm{\zeta}: ||\bm{\zeta}||_2 \leq 1
\end{aligned}
\end{align}
Now to make the solution $\bm x$ feasible over an ellipsoid of radius $\delta_E$, a transformation $\zeta_j \to \dfrac{\zeta_j}{\delta_E}$ is applied in \eqref{eq4.8}, which then becomes
\begin{align}
	\begin{aligned}	
		& \bm{x}^\top \left[\bm{\mu}^{(0)}\delta_{E}+\sum_{j=1}^{n}\zeta_j \bm{\mu}^{(j)}\right] \left[\bm{\mu}^{(0)}\delta_{E}+\sum_{j=1}^{n}\zeta_j \bm{\mu}^{(j)}\right]^\top\bm{x} - \tau^2 {\delta_E}^2 \geq 0, \qquad \forall \bm{\zeta}: ||\bm{\zeta}||_2 \leq \delta_{E}
	\end{aligned}
	\label{eq4.9}
\end{align}
Equivalently the above condition is written as,
\begin{align}
\begin{aligned}	\label{eq4.10}
& \bm{x}^\top \left[\bm{\mu}^{(0)}\delta_{E}+\sum_{j=1}^{n}\zeta_j \bm{\mu}^{(j)}\right] \left[\bm{\mu}^{(0)}\delta_{E}+\sum_{j=1}^{n}\zeta_j \bm{\mu}^{(j)}\right]^\top\bm{x} - \tau^2 \delta_{E}^2  \geq 0, \quad \forall \bm{\zeta}: ||\bm{\zeta}||_2^2 \leq \delta_E^2
\end{aligned}
\end{align}
The value of ${||\bm{\zeta}||_2}^2$ is $|{\zeta_1}|^2+|{\zeta_2}|^2+\dots+|{\zeta_n}|^2$ which can be written in a quadratic form as, $\bm{\zeta}^\top \bm{\zeta}$. So the condition \eqref{eq4.10} is equivalent to, 
\begin{align}
\begin{aligned}	\label{eq4.11}
& \bm{x}^\top \left[\bm{\mu}^{(0)}\delta_{E}+\sum_{j=1}^{n}\zeta_j \bm{\mu}^{(j)}\right] \left[\bm{\mu}^{(0)}\delta_{E}+\sum_{j=1}^{n}\zeta_j \bm{\mu}^{(j)}\right]^\top\bm{x} - \tau^2 \delta_{E}^2 \geq 0, \quad \forall \bm{\zeta}: \delta_{E}^2-\bm{\zeta}^\top \bm{\zeta} \geq 0
\end{aligned}
\end{align}
We claim that, the left hand side of the left expression in \eqref{eq4.11} is a homogeneous quadratic form of the vector $(\delta_E,\zeta_1, \zeta_2, \dots, \zeta_n)$. We prove this as follows.\\
To begin with, the vectors used in \eqref{eq4.10} need to be expressed componentwise for more clarity.\\
Now we have, \begin{align}
\begin{aligned}	\label{eq4.11a}
& \bm{\mu}^{(0)} \delta_E+\sum_{j=1}^{n}\zeta_j \bm{\mu}^{(j)} && = \begin{bmatrix} 
\mu_1^{(0)}\delta_E\\
\mu_2^{(0)}\delta_E\\
\vdots\\
\mu_n^{(0)}\delta_E
\end{bmatrix}
+ \sum_{j=1}^{n}\zeta_j  \begin{bmatrix} 
0\\
\vdots\\
\mu_j^{(j)}\\
\vdots\\
0
\end{bmatrix}\\
& && = \begin{bmatrix} 
\mu_1^{(0)}\delta_E+\zeta_1 \mu_1^{(1)}\\
\mu_2^{(0)}\delta_E+\zeta_2 \mu_1^{(2)}\\
\vdots\\
\mu_n^{(0)}\delta_E+\zeta_n \mu_n^{(n)}
\end{bmatrix}
\end{aligned}
\end{align}
Simplifying the left expression of \eqref{eq4.11} in component form, we get,
\begin{align}
\begin{aligned}	\label{eq4.11b}
\biggl [ \left(\mu_1^{(0)} \delta_E x_1+\mu_2^{(0)}\delta_E x_2+\dots+\mu_n^{(0)}\delta_E x_n\right)+\left(\zeta_1\mu_1^{(1)} x_1+\zeta_2\mu_2^{(2)} x_2+\dots+\zeta_n\mu_n^{(n)} x_n\right) \biggr ]^2-\tau^2 \delta_E^2\geq 0
\end{aligned}
\end{align}
Equivalently it can be written as,
\begin{align}
\begin{aligned}	\label{eq4.11c}
 &\left(\mu_1^{(0)} \delta_E x_1+\dots+\mu_n^{(0)}\delta_E x_n\right)^2-\tau^2\delta_E^2+\zeta_1^2\left(\mu_1^{(1)} x_1\right)^2+\dots+\zeta_n^2\left(\mu_n^{(n)} x_n\right)^2+2\sum_{i \neq j} \zeta_i \zeta_j \left(\mu_i^{(i)} x_i\right) \left(\mu_j^{(j)} x_j\right)\\
 & \qquad \qquad +2\left[\zeta_1 \mu_1^{(1)}x_1\left(\mu_1^{(0)}\delta_E x_1+\dots+\mu_n^{(0)}\delta_E x_n\right)+\dots+\zeta_n \mu_n^{(n)}x_n\left(\mu_1^{(0)}\delta_E x_1+\dots+\mu_n^{(0)}\delta_E x_n\right)\right]\geq 0
\end{aligned}
\end{align}
Now writing this in vector form, we get,
\begin{align}
\begin{aligned}	\label{eq4.11d}
&\delta_E^2\left[\left({\bm{\mu}^{(0)}}^\top\bm x\right)^2-\tau^2\right]+\zeta_1^2\left({\bm{\mu}^{(1)}}^\top\bm x\right)^2+\dots+\zeta_n^2\left({\bm{\mu}^{(n)}}^\top\bm x\right)^2+2\sum_{i \neq j} \zeta_i \zeta_j \left({\bm{\mu}^{(i)}}^\top\bm x\right)\left({\bm{\mu}^{(j)}}^\top\bm x\right)\\
& \qquad \qquad +2\left[\delta_E\zeta_1 \left({\bm{\mu}^{(1)}}^\top\bm x\right)\left({\bm{\mu}^{(0)}}^\top\bm x\right)+\dots+\delta_E\zeta_n \left({\bm{\mu}^{(n)}}^\top\bm x\right)\left({\bm{\mu}^{(0)}}^\top\bm x\right)\right]\geq 0
\end{aligned}
\end{align}
The above expression can be written as,
\begin{align}
\begin{aligned}	\label{eq4.11e}
\begin{scriptsize}
{\begin{bmatrix} 
	\delta_E\\
	\zeta_1\\
	\zeta_2\\
	\vdots\\
	\zeta_n
	\end{bmatrix}}^\top
\begin{bmatrix} 
({\bm{\mu}^{(0)}}^\top\bm x)^2-\tau^2 & ({\bm{\mu}^{(0)}}^\top\bm x)({\bm{\mu}^{(1)}}^\top\bm x)  & ({\bm{\mu}^{(0)}}^\top\bm x)({\bm{\mu}^{(2)}}^\top\bm x) & \dots & ({\bm{\mu}^{(0)}}^\top\bm x)({\bm{\mu}^{(n)}}^\top\bm x)\\
({\bm{\mu}^{(1)}}^\top\bm x)({\bm{\mu}^{(0)}}^\top\bm x) & ({\bm{\mu}^{(1)}}^\top\bm x)^2 & ({\bm{\mu}^{(1)}}^\top\bm x)({\bm{\mu}^{(2)}}^\top\bm x) & \dots & ({\bm{\mu}^{(1)}}^\top\bm x)({\bm{\mu}^{(n)}}^\top\bm x)\\
({\bm{\mu}^{(2)}}^\top\bm x)({\bm{\mu}^{(0)}}^\top\bm x) & ({\bm{\mu}^{(2)}}^\top\bm x)({\bm{\mu}^{(1)}}^\top\bm x) & ({\bm{\mu}^{(2)}}^\top\bm x)^2 & \dots & ({\bm{\mu}^{(2)}}^\top\bm x)({\bm{\mu}^{(n)}}^\top\bm x)\\
\vdots & \vdots & \vdots & \ddots & \vdots \\
({\bm{\mu}^{(n)}}^\top\bm x)({\bm{\mu}^{(0)}}^\top\bm x) & ({\bm{\mu}^{(n)}}^\top\bm x)({\bm{\mu}^{(1)}}^\top\bm x) &({\bm{\mu}^{(n)}}^\top\bm x)({\bm{\mu}^{(2)}}^\top\bm x) & \dots & ({\bm{\mu}^{(n)}}^\top\bm x)^2
\end{bmatrix}
\begin{bmatrix} 
\delta_E\\
\zeta_1\\
\zeta_2\\
\vdots\\
\zeta_n
\end{bmatrix} 
\end{scriptsize} \geq 0
\end{aligned}
\end{align}
Thus the left expression of \eqref{eq4.11} is given by,
\begin{align}
\begin{aligned}	\label{eq4.11f}
&(\bm z_E)^\top \bm A (\bm z_E) \geq 0,
\end{aligned}
\end{align}
where  \begin{align}
\begin{aligned} \label{eq4.11g}
\bm z_E=\begin{bmatrix} 
\delta_E\\
\zeta_1\\
\zeta_2\\
\vdots\\
\zeta_n
\end{bmatrix},\text{ and } \bm A=\begin{bmatrix} 
({\bm{\mu}^{(0)}}^\top\bm x)^2-\tau^2 & ({\bm{\mu}^{(0)}}^\top\bm x)({\bm{\mu}^{(1)}}^\top\bm x)  & ({\bm{\mu}^{(0)}}^\top\bm x)({\bm{\mu}^{(2)}}^\top\bm x) & \dots & ({\bm{\mu}^{(0)}}^\top\bm x)({\bm{\mu}^{(n)}}^\top\bm x)\\
({\bm{\mu}^{(1)}}^\top\bm x)({\bm{\mu}^{(0)}}^\top\bm x) & ({\bm{\mu}^{(1)}}^\top\bm x)^2 & ({\bm{\mu}^{(1)}}^\top\bm x)({\bm{\mu}^{(2)}}^\top\bm x) & \dots & ({\bm{\mu}^{(1)}}^\top\bm x)({\bm{\mu}^{(n)}}^\top\bm x)\\
({\bm{\mu}^{(2)}}^\top\bm x)({\bm{\mu}^{(0)}}^\top\bm x) & ({\bm{\mu}^{(2)}}^\top\bm x)({\bm{\mu}^{(1)}}^\top\bm x) & ({\bm{\mu}^{(2)}}^\top\bm x)^2 & \dots & ({\bm{\mu}^{(2)}}^\top\bm x)({\bm{\mu}^{(n)}}^\top\bm x)\\
\vdots & \vdots & \vdots & \ddots & \vdots \\
({\bm{\mu}^{(n)}}^\top\bm x)({\bm{\mu}^{(0)}}^\top\bm x) & ({\bm{\mu}^{(n)}}^\top\bm x)({\bm{\mu}^{(1)}}^\top\bm x) &({\bm{\mu}^{(n)}}^\top\bm x)({\bm{\mu}^{(2)}}^\top\bm x) & \dots & ({\bm{\mu}^{(n)}}^\top\bm x)^2
\end{bmatrix}
\end{aligned}
\end{align}
are respectively an $(n+1)$ dimensional vector and a symmetric matrix of order $(n+1) \times (n+1)$. \\
Moreover, the left hand side of \eqref{eq4.11d} is a quadratic form of $(\delta_E,\zeta_1, \zeta_2, \dots, \zeta_n)$, which proves our claim. Let us denote this quadratic form as, $Q_{\bm x}^{(E)}(\delta_E, \bm{\zeta})$. \\
Now writing the left hand side of the right expression of \eqref{eq4.11} in quadratic form, we have,
\begin{align}
\begin{aligned}	\label{eq4.11h}
& \delta_E^2-\bm{\zeta}^\top \bm{\zeta} &&= \delta_E^2-(\zeta_1^2+\zeta_2^2+\dots+\zeta_n^2)\\
& &&= {\begin{bmatrix} 
\delta_E\\
\zeta_1\\
\zeta_2\\
\vdots\\
\zeta_n
\end{bmatrix}}^\top
\begin{bmatrix} 
1 & 0  & 0 & \dots & 0\\
0 & -1 & 0 & \dots & 0\\
0 & 0 & -1 & \dots & 0\\
\vdots & \vdots & \vdots & \ddots & \vdots \\
0 & 0 &0 & \dots & -1
\end{bmatrix}
\begin{bmatrix} 
\delta_E\\
\zeta_1\\
\zeta_2\\
\vdots\\
\zeta_n
\end{bmatrix}
\end{aligned}
\end{align}
Let us denote this quadratic form as, $P^{(E)}(\delta_E, \bm{\zeta})$ and the right expression in \eqref{eq4.11} can be written as,
\begin{align}
\begin{aligned}	\label{eq4.11i}
& && {\begin{bmatrix} 
	\delta_E\\
	\zeta_1\\
	\zeta_2\\
	\vdots\\
	\zeta_n
	\end{bmatrix}}^\top
\begin{bmatrix} 
1 & 0  & 0 & \dots & 0\\
0 & -1 & 0 & \dots & 0\\
0 & 0 & -1 & \dots & 0\\
\vdots & \vdots & \vdots & \ddots & \vdots \\
0 & 0 &0 & \dots & -1
\end{bmatrix}
\begin{bmatrix} 
\delta_E\\
\zeta_1\\
\zeta_2\\
\vdots\\
\zeta_n
\end{bmatrix} \geq 0
\end{aligned}
\end{align}
This means the right expression of \eqref{eq4.11} is given by,
\begin{align}
\begin{aligned}	\label{eq4.11j}
&(\bm z_E)^\top \bm B^{(E)} (\bm z_E) \geq 0,
\end{aligned}
\end{align}
where 
\begin{align} 
\begin{aligned}	\label{eq4.11k}
\bm B^{(E)}=\begin{bmatrix} 
1 & 0  & 0 & \dots & 0\\
0 & -1 & 0 & \dots & 0\\
0 & 0 & -1 & \dots & 0\\
\vdots & \vdots & \vdots & \ddots & \vdots \\
0 & 0 &0 & \dots & -1
\end{bmatrix}
\end{aligned}
\end{align}
is a symmetric matrix of order $(n+1) \times (n+1)$.\\
Now using the expressions \eqref{eq4.11f} and \eqref{eq4.11j}, we can write \eqref{eq4.11} as,
 \begin{align}
 \begin{aligned}	\label{eq4.11l}
 & (\bm z_E)^\top \bm A (\bm z_E) \geq 0, \quad \forall (\bm z_E) : (\bm z_E)^\top \bm B^{(E)} (\bm z_E) \geq 0
 \end{aligned}
 \end{align}
 Then using Lemma \ref{lemma4.1} we have, the matrix $\bm A - \lambda \bm B^{(E)}$ is positive semidefinite.
 Thus the following matrix inequality holds.
\begin{align}
\begin{aligned} \label{eq4.11m}
\begin{bmatrix} 
({\bm{\mu}^{(0)}}^\top\bm x)^2-\tau^2-\lambda & ({\bm{\mu}^{(0)}}^\top\bm x)({\bm{\mu}^{(1)}}^\top\bm x)  & ({\bm{\mu}^{(0)}}^\top\bm x)({\bm{\mu}^{(2)}}^\top\bm x) & \dots & ({\bm{\mu}^{(0)}}^\top\bm x)({\bm{\mu}^{(n)}}^\top\bm x)\\
({\bm{\mu}^{(1)}}^\top\bm x)({\bm{\mu}^{(0)}}^\top\bm x) & ({\bm{\mu}^{(1)}}^\top\bm x)^2 +\lambda& ({\bm{\mu}^{(1)}}^\top\bm x)({\bm{\mu}^{(2)}}^\top\bm x) & \dots & ({\bm{\mu}^{(1)}}^\top\bm x)({\bm{\mu}^{(n)}}^\top\bm x)\\
({\bm{\mu}^{(2)}}^\top\bm x)({\bm{\mu}^{(0)}}^\top\bm x) & ({\bm{\mu}^{(2)}}^\top\bm x)({\bm{\mu}^{(1)}}^\top\bm x) & ({\bm{\mu}^{(2)}}^\top\bm x)^2+\lambda & \dots & ({\bm{\mu}^{(2)}}^\top\bm x)({\bm{\mu}^{(n)}}^\top\bm x)\\
\vdots & \vdots & \vdots & \ddots & \vdots \\
({\bm{\mu}^{(n)}}^\top\bm x)({\bm{\mu}^{(0)}}^\top\bm x) & ({\bm{\mu}^{(n)}}^\top\bm x)({\bm{\mu}^{(1)}}^\top\bm x) &({\bm{\mu}^{(n)}}^\top\bm x)({\bm{\mu}^{(2)}}^\top\bm x) & \dots & ({\bm{\mu}^{(n)}}^\top\bm x)^2+\lambda
\end{bmatrix} \succeq \bm 0
\end{aligned}
\end{align}  
Thus the feasibility condition of the robust counterpart problem \eqref{eq4.5} for ellipsoidal uncertainty can be expressed in the form of a positive semidefinite matrix. The result is given as the following theorem.
\begin{theorem} \label{thm4.1}
Consider the uncertain Markowitz portfolio problem \eqref{eq4.5} under ellipsoidal uncertainty. The solution vector $\bm x$ to the problem is feasible for the robust constraint in \eqref{eq4.5} if and only if there exists a $\lambda \geq 0$ such that the matrix 
\begin{scriptsize}
\begin{align*}
\begin{bmatrix} 
({\bm{\mu}^{(0)}}^\top\bm x)^2-\tau^2-\lambda & ({\bm{\mu}^{(0)}}^\top\bm x)({\bm{\mu}^{(1)}}^\top\bm x)  & ({\bm{\mu}^{(0)}}^\top\bm x)({\bm{\mu}^{(2)}}^\top\bm x) & \dots & ({\bm{\mu}^{(0)}}^\top\bm x)({\bm{\mu}^{(n)}}^\top\bm x)\\
({\bm{\mu}^{(1)}}^\top\bm x)({\bm{\mu}^{(0)}}^\top\bm x) & ({\bm{\mu}^{(1)}}^\top\bm x)^2+\lambda & ({\bm{\mu}^{(1)}}^\top\bm x)({\bm{\mu}^{(2)}}^\top\bm x) & \dots & ({\bm{\mu}^{(1)}}^\top\bm x)({\bm{\mu}^{(n)}}^\top\bm x)\\
({\bm{\mu}^{(2)}}^\top\bm x)({\bm{\mu}^{(0)}}^\top\bm x) & ({\bm{\mu}^{(2)}}^\top\bm x)({\bm{\mu}^{(1)}}^\top\bm x) & ({\bm{\mu}^{(2)}}^\top\bm x)^2+\lambda & \dots & ({\bm{\mu}^{(2)}}^\top\bm x)({\bm{\mu}^{(n)}}^\top\bm x)\\
\vdots & \vdots & \vdots & \ddots & \vdots \\
({\bm{\mu}^{(n)}}^\top\bm x)({\bm{\mu}^{(0)}}^\top\bm x) & ({\bm{\mu}^{(n)}}^\top\bm x)({\bm{\mu}^{(1)}}^\top\bm x) &({\bm{\mu}^{(n)}}^\top\bm x)({\bm{\mu}^{(2)}}^\top\bm x) & \dots & ({\bm{\mu}^{(n)}}^\top\bm x)^2+\lambda
\end{bmatrix}
\end{align*}
\end{scriptsize} is positive semidefinite.	
\end{theorem} 
\subsubsection{Box Uncertainty}
When the uncertain set in the problem \eqref{eq4.7} is a unit box, it is represented by,
\begin{align}
\begin{aligned}	\label{eq4.22}
\mathscr{U}_{\mu}^{(B)}=\left\{\bm{\mu} : \bm{\mu}= \bm{\mu}^{(0)}+\sum_{j=1}^{n}\zeta_j \bm{\mu}^{(j)},\ ||\bm{\zeta}||_\infty \leq 1\right\}
\end{aligned}
\end{align}
Then any solution $\bm{x}$ is robust feasible for the uncertain Markowitz potfolio model with box uncertainty, if and only if we have,
\begin{align}
\begin{aligned}	\label{eq4.23}
& \bm{x}^\top \left[\bm{\mu}^{(0)}+\sum_{j=1}^{n}\zeta_j \bm{\mu}^{(j)}\right] \left[\bm{\mu}^{(0)}+\sum_{j=1}^{n}\zeta_j \bm{\mu}^{(j)}\right]^\top\bm{x} - \tau^2 \geq 0, \quad \forall \bm{\zeta}: ||\bm{\zeta}||_{\infty} \leq 1
\end{aligned}
\end{align}
Now to make the solution $\bm x$ feasible over a box of radius $\delta_B$, a transformation $\zeta_j \to \dfrac{\zeta_j}{\delta_B}$ is applied in \eqref{eq4.23}, which then becomes
\begin{align}
\begin{aligned}	
& \bm{x}^\top \left[\bm{\mu}^{(0)}\delta_{B}+\sum_{j=1}^{n}\zeta_j \bm{\mu}^{(j)}\right] \left[\bm{\mu}^{(0)}\delta_{B}+\sum_{j=1}^{n}\zeta_j \bm{\mu}^{(j)} \right]^\top\bm{x} - \tau^2 {\delta_B}^2 \geq 0, \qquad \forall \bm{\zeta}: ||\bm{\zeta}||_{\infty} \leq \delta_{B}\\
\end{aligned}
\label{eq4.24}
\end{align}
Equivalently the above condition is written as,
\begin{align}
\begin{aligned}	\label{eq4.25}
& && \bm{x}^\top \left[\bm{\mu}^{(0)}\delta_{B}+\sum_{j=1}^{n}\zeta_j \bm{\mu}^{(j)}\right] \left[\bm{\mu}^{(0)}\delta_{B}+\sum_{j=1}^{n}\zeta_j \bm{\mu}^{(j)}\right]^\top\bm{x} - \tau^2 \delta_{B}^2 \geq 0, \quad \forall \bm{\zeta}:  {||\bm{\zeta}||_{\infty}}^2 \leq \delta_{B}^2\\
&\iff && \bm{x}^\top \left[\bm{\mu}^{(0)}\delta_{B}+\sum_{j=1}^{n}\zeta_j \bm{\mu}^{(j)}\right] \left[\bm{\mu}^{(0)}\delta_{B}+\sum_{j=1}^{n}\zeta_j \bm{\mu}^{(j)}\right]^\top\bm{x} - \tau^2 \delta_{B}^2  \geq 0, \quad \forall \bm{\zeta}:  \delta_{B}^2-{||\bm{\zeta}||_{\infty}}^2 \geq 0
\end{aligned}
\end{align}
Analogous to the ellipsoidal uncertainty case, the left hand side of the left expression of \eqref{eq4.25} can be written in quadratic form $Q_{\bm x}^{(B)}(\delta_{B}, \bm{\zeta})$, which is given by,
\begin{align}
\begin{aligned}	\label{eq4.25aa}
&(\bm z_B)^\top \bm A (\bm z_B) \geq 0,
\end{aligned}
\end{align}
where  \begin{align}
\begin{aligned} \label{eq4.25bb}
\bm z_B=\begin{bmatrix} 
\delta_B\\
\zeta_1\\
\zeta_2\\
\vdots\\
\zeta_n
\end{bmatrix},\text{ and } \bm A=\begin{bmatrix} 
({\bm{\mu}^{(0)}}^\top\bm x)^2-\tau^2 & ({\bm{\mu}^{(0)}}^\top\bm x)({\bm{\mu}^{(1)}}^\top\bm x)  & ({\bm{\mu}^{(0)}}^\top\bm x)({\bm{\mu}^{(2)}}^\top\bm x) & \dots & ({\bm{\mu}^{(0)}}^\top\bm x)({\bm{\mu}^{(n)}}^\top\bm x)\\
({\bm{\mu}^{(1)}}^\top\bm x)({\bm{\mu}^{(0)}}^\top\bm x) & ({\bm{\mu}^{(1)}}^\top\bm x)^2 & ({\bm{\mu}^{(1)}}^\top\bm x)({\bm{\mu}^{(2)}}^\top\bm x) & \dots & ({\bm{\mu}^{(1)}}^\top\bm x)({\bm{\mu}^{(n)}}^\top\bm x)\\
({\bm{\mu}^{(2)}}^\top\bm x)({\bm{\mu}^{(0)}}^\top\bm x) & ({\bm{\mu}^{(2)}}^\top\bm x)({\bm{\mu}^{(1)}}^\top\bm x) & ({\bm{\mu}^{(2)}}^\top\bm x)^2 & \dots & ({\bm{\mu}^{(2)}}^\top\bm x)({\bm{\mu}^{(n)}}^\top\bm x)\\
\vdots & \vdots & \vdots & \ddots & \vdots \\
({\bm{\mu}^{(n)}}^\top\bm x)({\bm{\mu}^{(0)}}^\top\bm x) & ({\bm{\mu}^{(n)}}^\top\bm x)({\bm{\mu}^{(1)}}^\top\bm x) &({\bm{\mu}^{(n)}}^\top\bm x)({\bm{\mu}^{(2)}}^\top\bm x) & \dots & ({\bm{\mu}^{(n)}}^\top\bm x)^2
\end{bmatrix}
\end{aligned}
\end{align}
are respectively an $(n+1)$ dimensional vector and a symmetric matrix of order $(n+1) \times (n+1)$.\\
Next the left hand side of the right expression of \eqref{eq4.25} needs to be written in quadratic form.
The value of ${||\bm{\zeta}||_{\infty}}^2$ equals to $\sup_j \zeta_j^2$. Let us denote it as ${\zeta_M}^2$. Now writing the left hand side of the right expression of \eqref{eq4.25} in quadratic form, we have,
\begin{align}
\begin{aligned}	\label{eq4.25a}
& \delta_{B}^2-{||\bm{\zeta}||_{\infty}}^2 &&= \delta_{B}^2-\zeta_M^2\\
& &&= {\begin{bmatrix} 
	\delta_{B}\\
	\zeta_1\\
	\vdots\\
	\zeta_M\\
	\vdots\\
	\zeta_n
	\end{bmatrix}}^\top
\begin{bmatrix} 
1 & \ 0& \ \dots  & \ 0& \ \dots & \ 0\\
0 & \ 0& \ \dots  & \ 0& \ \dots & \ 0\\
\vdots & \ \vdots & \ \ddots & \ \vdots & \ \ddots & \ \vdots \\
0 & \ 0 & \ \dots  & \ -1 & \ \dots & \ 0\\
\vdots & \ \vdots & \ \ddots & \ \vdots & \ \ddots & \ \vdots \\
0 & \ 0 & \dots  & \ 0 & \dots & \ 0
\end{bmatrix}
\begin{bmatrix} 
\delta_{B}\\
\zeta_1\\
\vdots\\
\zeta_M\\
\vdots\\
\zeta_n
\end{bmatrix}
\end{aligned}
\end{align}
Let us denote this quadratic form as, $P^{(B)}(\delta_{B}, \bm{\zeta})$. Now the right expression of \eqref{eq4.25} can be written as,
\begin{align}
\begin{aligned}	\label{eq4.25b}
&(\bm z_B)^\top \bm B^{(B)} (\bm z_B) \geq 0,
\end{aligned}
\end{align}
where 
\begin{align} 
\begin{aligned}	\label{eq4.25c}
\bm B^{(B)}=\begin{bmatrix} 
1 & \ 0& \ \dots  & \ 0& \ \dots & \ 0\\
0 & \ 0& \ \dots  & \ 0& \ \dots & \ 0\\
\vdots & \ \vdots & \ \ddots & \ \vdots & \ \ddots & \ \vdots \\
0 & \ 0 & \ \dots  & \ -1 & \ \dots & \ 0\\
\vdots & \ \vdots & \ \ddots & \ \vdots & \ \ddots & \ \vdots \\
0 & \ 0 & \dots  & \ 0 & \dots & \ 0
\end{bmatrix}
\end{aligned}
\end{align}
is a symmetric matrix of order $(n+1) \times (n+1)$. More specifically, $\bm {B^{(B)}}$ is a diagonal matrix with its first and $M^{th}$ diagonal elements are $1$ and $-1$ respectively, and the other diagonal elements are zeros.\\
Now using the expressions \eqref{eq4.25aa} and \eqref{eq4.25b}, we can write \eqref{eq4.25} as,
\begin{align}
	\begin{aligned}	\label{eq4.25d}
		& (\bm z_B)^\top \bm A \bm (z_B) \geq 0, \quad \forall (\bm z_B) : (\bm z_B)^\top \bm B^{(B)} (\bm z_B) \geq 0
	\end{aligned}
\end{align}
Then using Lemma \ref{lemma4.1} we have, the matrix $\bm A - \lambda \bm B^{(B)}$ is positive semidefinite.
Thus the following matrix inequality holds.
\begin{align} \label{eq4.25e}
%\begin{scriptsize}
\begin{bmatrix} 
	({\bm{\mu}^{(0)}}^\top\bm x)^2-\tau^2-\lambda & \ ({\bm{\mu}^{(0)}}^\top\bm x)({\bm{\mu}^{(1)}}^\top\bm x)& \ \dots  & \ ({\bm{\mu}^{(0)}}^\top\bm x)({\bm{\mu}^{(M)}}^\top\bm x)& \ \dots & \ ({\bm{\mu}^{(0)}}^\top\bm x)({\bm{\mu}^{(n)}}^\top\bm x)\\
	({\bm{\mu}^{(1)}}^\top\bm x)({\bm{\mu}^{(0)}}^\top\bm x) & \ ({\bm{\mu}^{(1)}}^\top\bm x)^2 & \ \dots & \ ({\bm{\mu}^{(1)}}^\top\bm x)({\bm{\mu}^{(M)}}^\top\bm x)& \ \dots  & \ ({\bm{\mu}^{(1)}}^\top\bm x)({\bm{\mu}^{(n)}}^\top\bm x)\\
	({\bm{\mu}^{(2)}}^\top\bm x)({\bm{\mu}^{(0)}}^\top\bm x) & \ ({\bm{\mu}^{(2)}}^\top\bm x)({\bm{\mu}^{(1)}}^\top\bm x)& \ \dots  & \ ({\bm{\mu}^{(2)}}^\top\bm x)({\bm{\mu}^{(M)}}^\top\bm x) & \ \dots & \ ({\bm{\mu}^{(2)}}^\top\bm x)({\bm{\mu}^{(n)}}^\top\bm x)\\
	\vdots & \ \vdots & \ \ddots & \ \vdots & \ \ddots & \ \vdots \\
	({\bm{\mu}^{(M)}}^\top\bm x)({\bm{\mu}^{(0)}}^\top\bm x) & \ ({\bm{\mu}^{(M)}}^\top\bm x)({\bm{\mu}^{(1)}}^\top\bm x)& \ \dots  & \ ({\bm{\mu}^{(M)}}^\top\bm x)^2+\lambda & \ \dots & \ ({\bm{\mu}^{(M)}}^\top\bm x)({\bm{\mu}^{(n)}}^\top\bm x)\\
	\vdots & \ \vdots & \ \ddots & \ \vdots & \ \ddots & \ \vdots \\
	({\bm{\mu}^{(n)}}^\top\bm x)({\bm{\mu}^{(0)}}^\top\bm x) & \ ({\bm{\mu}^{(n)}}^\top\bm x)({\bm{\mu}^{(1)}}^\top\bm x)& \ \dots  &\ ({\bm{\mu}^{(n)}}^\top\bm x)({\bm{\mu}^{(M)}}^\top\bm x) & \ \dots & \ ({\bm{\mu}^{(n)}}^\top\bm x)^2 
\end{bmatrix}
%\end{scriptsize}
\succeq \bm 0
\end{align} 
This can be written as the following result. 
\begin{theorem} \label{thm4.2}
Consider the uncertain portfolio problem \eqref{eq4.5} under box uncertainty. The solution vector $\bm x$ to the problem is feasible for the robust constraint in \eqref{eq4.5} if and only if there exists a $\lambda \geq 0$ such that the matrix 
\begin{scriptsize}
\begin{align*}
	\begin{bmatrix} 
		({\bm{\mu}^{(0)}}^\top\bm x)^2-\tau^2-\lambda & \ ({\bm{\mu}^{(0)}}^\top\bm x)({\bm{\mu}^{(1)}}^\top\bm x)& \ \dots  & \ ({\bm{\mu}^{(0)}}^\top\bm x)({\bm{\mu}^{(M)}}^\top\bm x)& \ \dots & \ ({\bm{\mu}^{(0)}}^\top\bm x)({\bm{\mu}^{(n)}}^\top\bm x)\\
		({\bm{\mu}^{(1)}}^\top\bm x)({\bm{\mu}^{(0)}}^\top\bm x) & \ ({\bm{\mu}^{(1)}}^\top\bm x)^2 & \ \dots & \ ({\bm{\mu}^{(1)}}^\top\bm x)({\bm{\mu}^{(M)}}^\top\bm x)& \ \dots  & \ ({\bm{\mu}^{(1)}}^\top\bm x)({\bm{\mu}^{(n)}}^\top\bm x)\\
		({\bm{\mu}^{(2)}}^\top\bm x)({\bm{\mu}^{(0)}}^\top\bm x) & \ ({\bm{\mu}^{(2)}}^\top\bm x)({\bm{\mu}^{(1)}}^\top\bm x)& \ \dots  & \ ({\bm{\mu}^{(2)}}^\top\bm x)({\bm{\mu}^{(M)}}^\top\bm x) & \ \dots & \ ({\bm{\mu}^{(2)}}^\top\bm x)({\bm{\mu}^{(n)}}^\top\bm x)\\
		\vdots & \ \vdots & \ \ddots & \ \vdots & \ \ddots & \ \vdots \\
		({\bm{\mu}^{(M)}}^\top\bm x)({\bm{\mu}^{(0)}}^\top\bm x) & \ ({\bm{\mu}^{(M)}}^\top\bm x)({\bm{\mu}^{(1)}}^\top\bm x)& \ \dots  & \ ({\bm{\mu}^{(M)}}^\top\bm x)^2+\lambda & \ \dots & \ ({\bm{\mu}^{(M)}}^\top\bm x)({\bm{\mu}^{(n)}}^\top\bm x)\\
		\vdots & \ \vdots & \ \ddots & \ \vdots & \ \ddots & \ \vdots \\
		({\bm{\mu}^{(n)}}^\top\bm x)({\bm{\mu}^{(0)}}^\top\bm x) & \ ({\bm{\mu}^{(n)}}^\top\bm x)({\bm{\mu}^{(1)}}^\top\bm x)& \ \dots  &\ ({\bm{\mu}^{(n)}}^\top\bm x)({\bm{\mu}^{(M)}}^\top\bm x) & \ \dots & \ ({\bm{\mu}^{(n)}}^\top\bm x)^2\\
	\end{bmatrix}
\end{align*}
\end{scriptsize} is positive semidefinite, where $M$ is the index of $\sup_j \zeta_j^2$.
\end{theorem}
\subsubsection{Polyhedral Uncertainty}
When the uncertain set in the problem \eqref{eq4.7} is a unit polyhedron, it is represented by,
\begin{align}
\begin{aligned}	\label{eq4.14}
\mathscr{U}_{\mu}^{(P)}=\left\{\bm{\mu} : \bm{\mu}= \bm{\mu}^{(0)}+\sum_{j=1}^{n}\zeta_j \bm{\mu}^{(j)},\ ||\bm{\zeta}||_1 \leq 1\right\}
\end{aligned}
\end{align}
Then any solution $\bm{x}$ is robust feasible for the uncertain Markowitz potfolio model with polyhedral uncertainty, if and only if we have,
\begin{align}
\begin{aligned}	\label{eq4.15}
& \bm{x}^\top \left[\bm{\mu}^{(0)}+\sum_{j=1}^{n}\zeta_j \bm{\mu}^{(j)}\right] \left[\bm{\mu}^{(0)}+\sum_{j=1}^{n}\zeta_j \bm{\mu}^{(j)}\right]^\top\bm{x} - \tau^2 \geq 0, \quad \forall \bm{\zeta}: ||\bm{\zeta}||_1 \leq 1
\end{aligned}
\end{align}
Now to make the solution $\bm x$ feasible over a polyhedron of radius $\delta_P$, a transformation $\zeta_j \to \dfrac{\zeta_j}{\delta_P}$ is applied in \eqref{eq4.23}, which then becomes
\begin{align}
\begin{aligned}	\label{eq4.15a}
	& \bm{x}^\top \left[\bm{\mu}^{(0)}\delta_{P}+\sum_{j=1}^{n}\zeta_j \bm{\mu}^{(j)}\right] \left[\bm{\mu}^{(0)}\delta_{P}+\sum_{j=1}^{n}\zeta_j \bm{\mu}^{(j)} \right]^\top\bm{x} - \tau^2 {\delta_P}^2 \geq 0, \qquad \forall \bm{\zeta}: ||\bm{\zeta}||_{1} \leq \delta_{P}\\
\end{aligned}
\end{align}
Equivalently the above condition is written as,
\begin{align}
\begin{aligned}	\label{eq4.16}
	& && \bm{x}^\top \left[\bm{\mu}^{(0)}\delta_{P}+\sum_{j=1}^{n}\zeta_j \bm{\mu}^{(j)}\right] \left[\bm{\mu}^{(0)}\delta_{P}+\sum_{j=1}^{n}\zeta_j \bm{\mu}^{(j)}\right]^\top\bm{x} - \tau^2 \delta_{P}^2  \geq 0, \quad \forall \bm{\zeta}:  {||\bm{\zeta}||_{1}}^2 \leq \delta_{P}^2\\
	&\iff && \bm{x}^\top \left[\bm{\mu}^{(0)}\delta_{P}+\sum_{j=1}^{n}\zeta_j \bm{\mu}^{(j)}\right] \left[\bm{\mu}^{(0)}\delta_{P}+\sum_{j=1}^{n}\zeta_j \bm{\mu}^{(j)}\right]^\top\bm{x} - \tau^2 \delta_{P}^2 \geq 0, \quad \forall \bm{\zeta}:  \delta_{P}^2-{||\bm{\zeta}||_{1}}^2 \geq 0
\end{aligned}
\end{align}
The left hand side of the left expression of \eqref{eq4.16} can be written in quadratic form $Q_{\bm x}^{(P)}(\delta_{P}, \bm{\zeta})$, which is given by,
\begin{align}
\begin{aligned}	\label{eq4.16a}
&(\bm z_P)^\top \bm A (\bm z_P) \geq 0,
\end{aligned}
\end{align}
where  \begin{align}
\begin{aligned} \label{eq4.16b}
\bm z_P=\begin{bmatrix} 
\delta_P\\
\zeta_1\\
\zeta_2\\
\vdots\\
\zeta_n
\end{bmatrix},\text{ and } \bm A=\begin{bmatrix} 
({\bm{\mu}^{(0)}}^\top\bm x)^2-\tau^2 & ({\bm{\mu}^{(0)}}^\top\bm x)({\bm{\mu}^{(1)}}^\top\bm x)  & ({\bm{\mu}^{(0)}}^\top\bm x)({\bm{\mu}^{(2)}}^\top\bm x) & \dots & ({\bm{\mu}^{(0)}}^\top\bm x)({\bm{\mu}^{(n)}}^\top\bm x)\\
({\bm{\mu}^{(1)}}^\top\bm x)({\bm{\mu}^{(0)}}^\top\bm x) & ({\bm{\mu}^{(1)}}^\top\bm x)^2 & ({\bm{\mu}^{(1)}}^\top\bm x)({\bm{\mu}^{(2)}}^\top\bm x) & \dots & ({\bm{\mu}^{(1)}}^\top\bm x)({\bm{\mu}^{(n)}}^\top\bm x)\\
({\bm{\mu}^{(2)}}^\top\bm x)({\bm{\mu}^{(0)}}^\top\bm x) & ({\bm{\mu}^{(2)}}^\top\bm x)({\bm{\mu}^{(1)}}^\top\bm x) & ({\bm{\mu}^{(2)}}^\top\bm x)^2 & \dots & ({\bm{\mu}^{(2)}}^\top\bm x)({\bm{\mu}^{(n)}}^\top\bm x)\\
\vdots & \vdots & \vdots & \ddots & \vdots \\
({\bm{\mu}^{(n)}}^\top\bm x)({\bm{\mu}^{(0)}}^\top\bm x) & ({\bm{\mu}^{(n)}}^\top\bm x)({\bm{\mu}^{(1)}}^\top\bm x) &({\bm{\mu}^{(n)}}^\top\bm x)({\bm{\mu}^{(2)}}^\top\bm x) & \dots & ({\bm{\mu}^{(n)}}^\top\bm x)^2
\end{bmatrix}
\end{aligned}
\end{align}
are respectively an $(n+1)$ dimensional vector and a symmetric matrix of order $(n+1) \times (n+1)$.\\
The term ${||\bm{\zeta}||_1}^2$ in the right expression of \eqref{eq4.16} can be written as a quadratic form as $\bm{\zeta}^\top S  \bm{\zeta}$, where $S$ is an $(n \times n)$ symmetric matrix with $1$ as each entry and is given by,
\begin{align}
\begin{aligned} \label{eq4.16c}
S= 	\begin{bmatrix} 
1 && 1  && 1 && \dots && 1\\
1 && 1  && 1 && \dots && 1\\
\vdots && \vdots && \vdots && \ddots && \vdots \\
1 && 1 && 1 && \dots && 1\\
\end{bmatrix}
\end{aligned} 
\end{align}
and hence the left hand side of the right expression of \eqref{eq4.16} can be written in quadratic form as,
\begin{align}
\begin{aligned}	\label{eq4.17}
& \delta_{P}^2-{||\bm{\zeta}||_{1}}^2 &&= {\begin{bmatrix} 
	\delta_{P}\\
	\zeta_1\\
	\zeta_2\\
	\vdots\\
	\zeta_n
\end{bmatrix}}^\top
\begin{bmatrix} 
	1 & 0  & 0 & \dots & 0\\
	0 & -1 & -1 & \dots & -1\\
	0 & -1 & -1 & \dots & -1\\
	\vdots & \vdots & \vdots & \ddots & \vdots \\
	0 & -1 & -1 & \dots & -1
\end{bmatrix}
\begin{bmatrix} 
	\delta_{P}\\
	\zeta_1\\
	\zeta_2\\
	\vdots\\
	\zeta_n
\end{bmatrix}
\end{aligned}
\end{align}
Let us denote this quadratic form as, $P^{(P)}(\delta_{P}, \bm{\zeta})$. Now the right expression of \eqref{eq4.16} can be written as,
\begin{align}
\begin{aligned}	\label{eq4.18}
&(\bm z_P)^\top \bm B^{(P)} \bm (z_P) \geq 0,
\end{aligned}
\end{align}
where 
\begin{align} 
\begin{aligned}	\label{eq4.19}
\bm B^{(P)}=\begin{bmatrix} 
	1 & 0  & 0 & \dots & 0\\
	0 & -1 & -1 & \dots & -1\\
	0 & -1 & -1 & \dots & -1\\
	\vdots & \vdots & \vdots & \ddots & \vdots \\
	0 & -1 & -1 & \dots & -1
\end{bmatrix}
\end{aligned}
\end{align}
is a symmetric matrix of order $(n+1) \times (n+1)$.\\
Moreover, the left expression of \eqref{eq4.16} is given by,
\begin{align}
\begin{aligned}	\label{eq4.20}
&(\bm z_P)^\top \bm A (\bm z_P) \geq 0,
\end{aligned}
\end{align}
where $\bm A$ is the $(n+1) \times (n+1)$ symmetric matrix given in \eqref{eq4.11g}.\\
Now using the expressions \eqref{eq4.18} and \eqref{eq4.20}, we can write \eqref{eq4.16} as,
\begin{align}
\begin{aligned}	\label{eq4.21}
& (\bm z_P)^\top \bm A (\bm z_P) \geq 0, \quad \forall \bm z_P : (\bm z_P)^\top \bm B^{(P)} (\bm z_P) \geq 0
\end{aligned}
\end{align}
Then applying Lemma \ref{lemma4.1}, we have the matrix $\bm A - \lambda \bm B^{(P)}$ is positive semidefinite and hence the following matrix inequality holds. 
\begin{align*}
\begin{bmatrix} 
({\bm{\mu}^{(0)}}^\top\bm x)^2-\tau^2-\lambda & ({\bm{\mu}^{(0)}}^\top\bm x)({\bm{\mu}^{(1)}}^\top\bm x)  & ({\bm{\mu}^{(0)}}^\top\bm x)({\bm{\mu}^{(2)}}^\top\bm x) & \dots & ({\bm{\mu}^{(0)}}^\top\bm x)({\bm{\mu}^{(n)}}^\top\bm x)\\
({\bm{\mu}^{(1)}}^\top\bm x)({\bm{\mu}^{(0)}}^\top\bm x) & ({\bm{\mu}^{(1)}}^\top\bm x)^2+\lambda & ({\bm{\mu}^{(1)}}^\top\bm x)({\bm{\mu}^{(2)}}^\top\bm x)+\lambda & \dots & ({\bm{\mu}^{(1)}}^\top\bm x)({\bm{\mu}^{(n)}}^\top\bm x)+\lambda\\
({\bm{\mu}^{(2)}}^\top\bm x)({\bm{\mu}^{(0)}}^\top\bm x) & ({\bm{\mu}^{(2)}}^\top\bm x)({\bm{\mu}^{(1)}}^\top\bm x)+\lambda & ({\bm{\mu}^{(2)}}^\top\bm x)^2+\lambda & \dots & ({\bm{\mu}^{(2)}}^\top\bm x)({\bm{\mu}^{(n)}}^\top\bm x)+\lambda\\
\vdots & \vdots & \vdots & \ddots & \vdots \\
({\bm{\mu}^{(n)}}^\top\bm x)({\bm{\mu}^{(0)}}^\top\bm x) & ({\bm{\mu}^{(n)}}^\top\bm x)({\bm{\mu}^{(1)}}^\top\bm x)+\lambda &({\bm{\mu}^{(n)}}^\top\bm x)({\bm{\mu}^{(2)}}^\top\bm x)+\lambda & \dots & ({\bm{\mu}^{(n)}}^\top\bm x)^2+\lambda 
\end{bmatrix} \succeq \bm 0
\end{align*}
The above result is given as the following theorem.
\begin{theorem} \label{thm4.3}
Consider the uncertain portfolio problem \eqref{eq4.5} under polyhedral uncertainty. The solution vector $\bm x$ to the problem is feasible for the robust constraint in \eqref{eq4.5} if and only if there exists a $\lambda \geq 0$ such that the matrix 
\begin{scriptsize}
\begin{align*}
\begin{bmatrix} 
({\bm{\mu}^{(0)}}^\top\bm x)^2-\tau^2-\lambda & ({\bm{\mu}^{(0)}}^\top\bm x)({\bm{\mu}^{(1)}}^\top\bm x)  & ({\bm{\mu}^{(0)}}^\top\bm x)({\bm{\mu}^{(2)}}^\top\bm x) & \dots & ({\bm{\mu}^{(0)}}^\top\bm x)({\bm{\mu}^{(n)}}^\top\bm x)\\
({\bm{\mu}^{(1)}}^\top\bm x)({\bm{\mu}^{(0)}}^\top\bm x) & ({\bm{\mu}^{(1)}}^\top\bm x)^2+\lambda & ({\bm{\mu}^{(1)}}^\top\bm x)({\bm{\mu}^{(2)}}^\top\bm x)+\lambda & \dots & ({\bm{\mu}^{(1)}}^\top\bm x)({\bm{\mu}^{(n)}}^\top\bm x)+\lambda\\
({\bm{\mu}^{(2)}}^\top\bm x)({\bm{\mu}^{(0)}}^\top\bm x) & ({\bm{\mu}^{(2)}}^\top\bm x)({\bm{\mu}^{(1)}}^\top\bm x)+\lambda & ({\bm{\mu}^{(2)}}^\top\bm x)^2+\lambda & \dots & ({\bm{\mu}^{(2)}}^\top\bm x)({\bm{\mu}^{(n)}}^\top\bm x)+\lambda\\
\vdots & \vdots & \vdots & \ddots & \vdots \\
({\bm{\mu}^{(n)}}^\top\bm x)({\bm{\mu}^{(0)}}^\top\bm x) & ({\bm{\mu}^{(n)}}^\top\bm x)({\bm{\mu}^{(1)}}^\top\bm x)+\lambda &({\bm{\mu}^{(n)}}^\top\bm x)({\bm{\mu}^{(2)}}^\top\bm x)+\lambda & \dots & ({\bm{\mu}^{(n)}}^\top\bm x)^2+\lambda
\end{bmatrix}
\end{align*}
\end{scriptsize} is positive semidefinite.	
\end{theorem}
\subsection{Feasibility Conditions for the Combined Uncertainties}
Being the intersection of two convex sets, each combined uncertain set is itself a convex set and can be represented as a quadratic form. Now we derive the feasibility conditions for the combined uncertainties.
\subsubsection{"Box $\cap$ Ellipsoidal" Uncertainty}
When the uncertain set in the problem \eqref{eq4.7} is the intersection of a unit box and a unit ellipsoid, it is represented by,
\begin{align}
\begin{aligned}	\label{eq4.37}
\mathscr{U}_{\mu}=\left\{\bm{\mu} : \bm{\mu}= \bm{\mu}^{(0)}+\sum_{j=1}^{n}\zeta_j \bm{\mu}^{(j)},\ ||\bm{\zeta}||_{\infty} \leq 1, \ ||\bm{\zeta}||_2 \leq 1 \right\}
\end{aligned}
\end{align}
Then any solution $\bm{x}$ is robust feasible for the uncertain Markowitz potfolio model with "box $\cap$ ellipsoidal" uncertainty, if and only if we have,
\begin{align}
\begin{aligned}	\label{eq4.38}
& \bm{x}^\top \left[\bm{\mu}^{(0)}+\sum_{j=1}^{n}\zeta_j \bm{\mu}^{(j)}\right] \left[\bm{\mu}^{(0)}+\sum_{j=1}^{n}\zeta_j \bm{\mu}^{(j)} \right]^\top\bm{x} - \tau^2 \geq 0, \quad \forall \bm{\zeta}: ||\bm{\zeta}||_{\infty} \leq 1, \ ||\bm{\zeta}||_2 \leq 1
\end{aligned}
\end{align}
Before proceeding further, we need to find the intersection between the box and ellipsoidal sets of radii $\delta_B$ and $\delta_E$ respectively.\\
Now the box uncertain set of radius $\delta_B$ is given by,
\begin{align*}
\left\{\bm{\zeta} : ||\bm{\zeta}||_{\infty} \leq \delta_B \right\} \ \implies \left\{\zeta_j:|\zeta_j|\leq \delta_B, \quad \forall j\right\}
\end{align*}
Further, the ellipsoidal uncertain set of radius $\delta_E$ is given by,
\begin{align*}
\left\{\bm{\zeta} : ||\bm{\zeta}||_{2} \leq \delta_E\right\} \ \implies \left\{\zeta_j:\sum_{j=1}^{n}(\zeta_j)^2 \leq (\delta_E)^2 \right\}
\end{align*}
So the intersection of the two sets is written as,
\begin{align*}
\begin{aligned}
& &&\left\{\zeta_j:\sum_{j=1}^{n}(\zeta_j)^2 \leq \min\left\{n\delta_B^2, \delta_E^2\right\}\right\} \\
& \text{ or, } && \left\{ \bm{\zeta}: ||\bm{\zeta}||_{2}^2 \leq \min\left\{n\delta_B^2, \delta_E^2\right\} \right\}
 \end{aligned} 
\end{align*}
Now there arises two cases.\\
\underline{\textit{\textbf{Case I:} When $\min\left\{n\delta_B^2, \delta_E^2\right\}=n\delta_B^2$}}\\
For this case the intersection set satisfies $||\bm{\zeta}||_{2} \leq \sqrt n\delta_B$, and it can be achieved by making a transformation $\zeta_j \to \dfrac{\zeta_j}{\sqrt n \delta_B}$ in the constraint \eqref{eq4.38}, which then becomes,
\begin{align}
\begin{aligned}	\label{eq4.39}
& \bm{x}^\top \left[\bm{\mu}^{(0)}(\sqrt n\delta_B)+\sum_{j=1}^{n}\zeta_j \bm{\mu}^{(j)}\right] \left[\bm{\mu}^{(0)}(\sqrt n\delta_B)+\sum_{j=1}^{n}\zeta_j \bm{\mu}^{(j)}\right]^\top\bm{x} - \tau^2 n\delta_B^2 \geq 0, \quad \forall  \bm{\zeta}: ||\bm{\zeta}||_2^2 \leq n\delta_B^2
\end{aligned}
\end{align}
Equivalently the above condition is written as,
\begin{align}
\begin{aligned} \label{eq4.40}
\bm{x}^\top \left[\bm{\mu}^{(0)}(\sqrt n\delta_B)+\sum_{j=1}^{n}\zeta_j \bm{\mu}^{(j)}\right] \left[\bm{\mu}^{(0)}(\sqrt n\delta_B)+\sum_{j=1}^{n}\zeta_j \bm{\mu}^{(j)}\right]^\top\bm{x} - \tau^2 n\delta_B^2 \geq 0, \quad \forall \bm{\zeta}: n\delta_B^2-\bm{\zeta}^\top \bm{\zeta} \geq 0
\end{aligned}
\end{align}
Similar derivations as used in the single uncertainty cases can convert the left hand side of the left expression of \eqref{eq4.25} into quadratic form $Q_{\bm x}^{(BE)}(\delta_{B}, \bm{\zeta})$, which is given by,
\begin{align}
\begin{aligned}	\label{eq4.40a}
&(\bm z_B)^\top \bm A' (\bm z_B) \geq 0,
\end{aligned}
\end{align}
where  \begin{align}
\begin{aligned} \label{eq4.40b}
&\bm z_B=\begin{bmatrix} 
\delta_B\\
\zeta_1\\
\zeta_2\\
\vdots\\
\zeta_n
\end{bmatrix}\text{ is an $(n+1)$ dimensional vector, and }\\
&\bm A'=\begin{bmatrix} 
(\sqrt n{\bm{\mu}^{(0)}}^\top\bm x)^2-n \tau^2 & (\sqrt n{\bm{\mu}^{(0)}}^\top\bm x)({\bm{\mu}^{(1)}}^\top\bm x)  & (\sqrt n{\bm{\mu}^{(0)}}^\top\bm x)({\bm{\mu}^{(2)}}^\top\bm x) & \dots & (\sqrt n{\bm{\mu}^{(0)}}^\top\bm x)({\bm{\mu}^{(n)}}^\top\bm x)\\
({\bm{\mu}^{(1)}}^\top\bm x)(\sqrt n{\bm{\mu}^{(0)}}^\top\bm x) & ({\bm{\mu}^{(1)}}^\top\bm x)^2 & ({\bm{\mu}^{(1)}}^\top\bm x)({\bm{\mu}^{(2)}}^\top\bm x) & \dots & ({\bm{\mu}^{(1)}}^\top\bm x)({\bm{\mu}^{(n)}}^\top\bm x)\\
({\bm{\mu}^{(2)}}^\top\bm x)(\sqrt n{\bm{\mu}^{(0)}}^\top\bm x) & ({\bm{\mu}^{(2)}}^\top\bm x)({\bm{\mu}^{(1)}}^\top\bm x) & ({\bm{\mu}^{(2)}}^\top\bm x)^2 & \dots & ({\bm{\mu}^{(2)}}^\top\bm x)({\bm{\mu}^{(n)}}^\top\bm x)\\
\vdots & \vdots & \vdots & \ddots & \vdots \\
({\bm{\mu}^{(n)}}^\top\bm x)(\sqrt n{\bm{\mu}^{(0)}}^\top\bm x) & ({\bm{\mu}^{(n)}}^\top\bm x)({\bm{\mu}^{(1)}}^\top\bm x) &({\bm{\mu}^{(n)}}^\top\bm x)({\bm{\mu}^{(2)}}^\top\bm x) & \dots & ({\bm{\mu}^{(n)}}^\top\bm x)^2
\end{bmatrix}
\end{aligned}
\end{align}
is a symmetric matrix of order $(n+1) \times (n+1)$.\\
Now writing the left hand side of the right expression of \eqref{eq4.40} in quadratic form, we have,
\begin{align}
\begin{aligned}	\label{eq4.40c}
& n\delta_B^2-\bm{\zeta}^\top \bm{\zeta} &&= n\delta_B^2-(\zeta_1^2+\zeta_2^2+\dots+\zeta_n^2)\\
& &&= {\begin{bmatrix} 
\delta_B\\
\zeta_1\\
\zeta_2\\
\vdots\\
\zeta_n
\end{bmatrix}}^\top
\begin{bmatrix} 
n & 0  & 0 & \dots & 0\\
0 & -1 & 0 & \dots & 0\\
0 & 0 & -1 & \dots & 0\\
\vdots & \vdots & \vdots & \ddots & \vdots \\
0 & 0 &0 & \dots & -1
\end{bmatrix}
\begin{bmatrix} 
\delta_B\\
\zeta_1\\
\zeta_2\\
\vdots\\
\zeta_n
\end{bmatrix}
\end{aligned}
\end{align}
Let us denote this quadratic form as, $P_1^{(BE)}(\delta_B, \bm{\zeta})$ and the right expression in \eqref{eq4.40} can be written as,
\begin{align}
\begin{aligned}	\label{eq4.41}
& && {\begin{bmatrix} 
\delta_B\\
\zeta_1\\
\zeta_2\\
\vdots\\
\zeta_n
\end{bmatrix}}^\top
\begin{bmatrix} 
n & 0  & 0 & \dots & 0\\
0 & -1 & 0 & \dots & 0\\
0 & 0 & -1 & \dots & 0\\
\vdots & \vdots & \vdots & \ddots & \vdots \\
0 & 0 &0 & \dots & -1
\end{bmatrix}
\begin{bmatrix} 
\delta_B\\
\zeta_1\\
\zeta_2\\
\vdots\\
\zeta_n
\end{bmatrix} \geq 0
\end{aligned}
\end{align}
That means,
\begin{align}
\begin{aligned}	\label{eq4.41a}
&(\bm z_B)^\top \bm B_1^{(BE)} (\bm z_B) \geq 0,
\end{aligned}
\end{align}
where 
\begin{align} 
\begin{aligned}	\label{eq4.41b}
\bm B_1^{(BE)}=\begin{bmatrix} 
n & 0  & 0 & \dots & 0\\
0 & -1 & 0 & \dots & 0\\
0 & 0 & -1 & \dots & 0\\
\vdots & \vdots & \vdots & \ddots & \vdots \\
0 & 0 &0 & \dots & -1
\end{bmatrix}
\end{aligned}
\end{align}
is a symmetric matrix of order $(n+1) \times (n+1)$.\\
Now using the expressions \eqref{eq4.40a} and \eqref{eq4.41a}, we can write \eqref{eq4.40} as,
\begin{align}
\begin{aligned}	\label{eq4.41c}
& (\bm z_B)^\top \bm A' (\bm z_B) \geq 0, \quad \forall (\bm z_B) : (\bm z_B)^\top \bm B_1^{(BE)} (\bm z_B) \geq 0
\end{aligned}
\end{align}
Then by Lemma \ref{lemma4.1} we have, the matrix $\bm A' - \lambda \bm B_1^{(BE)}$ is positive semidefinite.
Thus the following matrix inequality holds.
Therefore the matrix
\begin{align}
\begin{aligned}	\label{eq4.41d}
\begin{scriptsize}
\begin{bmatrix} 
(\sqrt n{\bm{\mu}^{(0)}}^\top\bm x)^2-n \tau^2-n\lambda & (\sqrt n{\bm{\mu}^{(0)}}^\top\bm x)({\bm{\mu}^{(1)}}^\top\bm x)  & (\sqrt n{\bm{\mu}^{(0)}}^\top\bm x)({\bm{\mu}^{(2)}}^\top\bm x) & \dots & (\sqrt n{\bm{\mu}^{(0)}}^\top\bm x)({\bm{\mu}^{(n)}}^\top\bm x)\\
({\bm{\mu}^{(1)}}^\top\bm x)(\sqrt n{\bm{\mu}^{(0)}}^\top\bm x) & ({\bm{\mu}^{(1)}}^\top\bm x)^2+\lambda & ({\bm{\mu}^{(1)}}^\top\bm x)({\bm{\mu}^{(2)}}^\top\bm x) & \dots & ({\bm{\mu}^{(1)}}^\top\bm x)({\bm{\mu}^{(n)}}^\top\bm x)\\
({\bm{\mu}^{(2)}}^\top\bm x)(\sqrt n{\bm{\mu}^{(0)}}^\top\bm x) & ({\bm{\mu}^{(2)}}^\top\bm x)({\bm{\mu}^{(1)}}^\top\bm x) & ({\bm{\mu}^{(2)}}^\top\bm x)^2+\lambda & \dots & ({\bm{\mu}^{(2)}}^\top\bm x)({\bm{\mu}^{(n)}}^\top\bm x)\\
\vdots & \vdots & \vdots & \ddots & \vdots \\
({\bm{\mu}^{(n)}}^\top\bm x)(\sqrt n{\bm{\mu}^{(0)}}^\top\bm x) & ({\bm{\mu}^{(n)}}^\top\bm x)({\bm{\mu}^{(1)}}^\top\bm x) &({\bm{\mu}^{(n)}}^\top\bm x)({\bm{\mu}^{(2)}}^\top\bm x) & \dots & ({\bm{\mu}^{(n)}}^\top\bm x)^2+\lambda
\end{bmatrix}
\end{scriptsize} \succeq \bm 0
\end{aligned}
\end{align} 
\\
\underline{\textit{\textbf{Case II:} When $\min\left\{n\delta_B^2, \delta_E^2\right\}=\delta_E^2$}}\\
For this case the intersection set satisfies $||\bm{\zeta}||_{2} \leq \delta_E$, and it can be achieved by making a transformation $\zeta_j \to \dfrac{\zeta_j}{\delta_E}$ in the constraint \eqref{eq4.38}, which then becomes,
\begin{align}
\begin{aligned}	\label{eq4.42}
& \bm{x}^\top \left[\bm{\mu}^{(0)}(\delta_E)+\sum_{j=1}^{n}\zeta_j \bm{\mu}^{(j)}\right] \left[\bm{\mu}^{(0)}(\delta_E)+\sum_{j=1}^{n}\zeta_j \bm{\mu}^{(j)}\right]^\top\bm{x} - \tau^2 \delta_E^2 \geq 0, \qquad \forall \bm{\zeta}: ||\bm{\zeta}||_2^2 \leq \delta_E^2
\end{aligned}
\end{align}
This is equivalent to the ellipsoidal uncertainty case given in Theorem \ref{thm4.1}. So proceeding in a similar manner, the following matrix inequality holds.
\begin{align}
\begin{aligned}	\label{eq4.42a}
\begin{scriptsize}
\begin{bmatrix} 
({\bm{\mu}^{(0)}}^\top\bm x)^2-\tau^2-\lambda & ({\bm{\mu}^{(0)}}^\top\bm x)({\bm{\mu}^{(1)}}^\top\bm x)  & ({\bm{\mu}^{(0)}}^\top\bm x)({\bm{\mu}^{(2)}}^\top\bm x) & \dots & ({\bm{\mu}^{(0)}}^\top\bm x)({\bm{\mu}^{(n)}}^\top\bm x)\\
({\bm{\mu}^{(1)}}^\top\bm x)({\bm{\mu}^{(0)}}^\top\bm x) & ({\bm{\mu}^{(1)}}^\top\bm x)^2+\lambda & ({\bm{\mu}^{(1)}}^\top\bm x)({\bm{\mu}^{(2)}}^\top\bm x) & \dots & ({\bm{\mu}^{(1)}}^\top\bm x)({\bm{\mu}^{(n)}}^\top\bm x)\\
({\bm{\mu}^{(2)}}^\top\bm x)({\bm{\mu}^{(0)}}^\top\bm x) & ({\bm{\mu}^{(2)}}^\top\bm x)({\bm{\mu}^{(1)}}^\top\bm x)& ({\bm{\mu}^{(2)}}^\top\bm x)^2+\lambda & \dots & ({\bm{\mu}^{(2)}}^\top\bm x)({\bm{\mu}^{(n)}}^\top\bm x)\\
\vdots & \vdots & \vdots & \ddots & \vdots \\
({\bm{\mu}^{(n)}}^\top\bm x)({\bm{\mu}^{(0)}}^\top\bm x) & ({\bm{\mu}^{(n)}}^\top\bm x)({\bm{\mu}^{(1)}}^\top\bm x) &({\bm{\mu}^{(n)}}^\top\bm x)({\bm{\mu}^{(2)}}^\top\bm x) & \dots & ({\bm{\mu}^{(n)}}^\top\bm x)^2+\lambda
\end{bmatrix}
\end{scriptsize} \succeq \bm 0
\end{aligned}
\end{align}
The results of the above two cases can be given as follows.
\begin{theorem} \label{thm4.4}
Consider the uncertain portfolio problem \eqref{eq4.5} under the uncertainty set obtained from the combination of box and ellipsoidal uncertainty sets. Let the radii of the box and ellipsoidal sets are respectively $\delta_B$ and $\delta_E$. Then the solution vector $\bm x$ to this problem is feasible for the robust constraint in \eqref{eq4.5} if and only if there exists a $\lambda \geq 0$ such that the following matrix inequalities hold:
\begin{tiny}
\begin{align*}
\begin{aligned}	
& \left \{ \begin{array}{ll}
\begin{bmatrix} 
(\sqrt n{\bm{\mu}^{(0)}}^\top\bm x)^2-n \tau^2-n\lambda & (\sqrt n{\bm{\mu}^{(0)}}^\top\bm x)({\bm{\mu}^{(1)}}^\top\bm x)  & (\sqrt n{\bm{\mu}^{(0)}}^\top\bm x)({\bm{\mu}^{(2)}}^\top\bm x) & \dots & (\sqrt n{\bm{\mu}^{(0)}}^\top\bm x)({\bm{\mu}^{(n)}}^\top\bm x)\\
({\bm{\mu}^{(1)}}^\top\bm x)(\sqrt n{\bm{\mu}^{(0)}}^\top\bm x) & ({\bm{\mu}^{(1)}}^\top\bm x)^2+\lambda & ({\bm{\mu}^{(1)}}^\top\bm x)({\bm{\mu}^{(2)}}^\top\bm x) & \dots & ({\bm{\mu}^{(1)}}^\top\bm x)({\bm{\mu}^{(n)}}^\top\bm x)\\
({\bm{\mu}^{(2)}}^\top\bm x)(\sqrt n{\bm{\mu}^{(0)}}^\top\bm x) & ({\bm{\mu}^{(2)}}^\top\bm x)({\bm{\mu}^{(1)}}^\top\bm x) & ({\bm{\mu}^{(2)}}^\top\bm x)^2+\lambda & \dots & ({\bm{\mu}^{(2)}}^\top\bm x)({\bm{\mu}^{(n)}}^\top\bm x)\\
\vdots & \vdots & \vdots & \ddots & \vdots \\
({\bm{\mu}^{(n)}}^\top\bm x)(\sqrt n{\bm{\mu}^{(0)}}^\top\bm x) & ({\bm{\mu}^{(n)}}^\top\bm x)({\bm{\mu}^{(1)}}^\top\bm x) &({\bm{\mu}^{(n)}}^\top\bm x)({\bm{\mu}^{(2)}}^\top\bm x) & \dots & ({\bm{\mu}^{(n)}}^\top\bm x)^2+\lambda\\
\end{bmatrix} \succeq 0, &\mbox{ if } \min\left\{n\delta_B^2, \delta_E^2\right\}=n\delta_B^2\\
\\
\begin{bmatrix} 
({\bm{\mu}^{(0)}}^\top\bm x)^2-\tau^2-\lambda & ({\bm{\mu}^{(0)}}^\top\bm x)({\bm{\mu}^{(1)}}^\top\bm x)  & ({\bm{\mu}^{(0)}}^\top\bm x)({\bm{\mu}^{(2)}}^\top\bm x) & \dots & ({\bm{\mu}^{(0)}}^\top\bm x)({\bm{\mu}^{(n)}}^\top\bm x)\\
({\bm{\mu}^{(1)}}^\top\bm x)({\bm{\mu}^{(0)}}^\top\bm x) & ({\bm{\mu}^{(1)}}^\top\bm x)^2+\lambda & ({\bm{\mu}^{(1)}}^\top\bm x)({\bm{\mu}^{(2)}}^\top\bm x) & \dots & ({\bm{\mu}^{(1)}}^\top\bm x)({\bm{\mu}^{(n)}}^\top\bm x)\\
({\bm{\mu}^{(2)}}^\top\bm x)({\bm{\mu}^{(0)}}^\top\bm x) & ({\bm{\mu}^{(2)}}^\top\bm x)({\bm{\mu}^{(1)}}^\top\bm x)& ({\bm{\mu}^{(2)}}^\top\bm x)^2+\lambda & \dots & ({\bm{\mu}^{(2)}}^\top\bm x)({\bm{\mu}^{(n)}}^\top\bm x)\\
\vdots & \vdots & \vdots & \ddots & \vdots \\
({\bm{\mu}^{(n)}}^\top\bm x)({\bm{\mu}^{(0)}}^\top\bm x) & ({\bm{\mu}^{(n)}}^\top\bm x)({\bm{\mu}^{(1)}}^\top\bm x) &({\bm{\mu}^{(n)}}^\top\bm x)({\bm{\mu}^{(2)}}^\top\bm x) & \dots & ({\bm{\mu}^{(n)}}^\top\bm x)^2+\lambda
\end{bmatrix} \succeq 0, &\mbox{ if } \min\left\{n\delta_B^2, \delta_E^2\right\}=\delta_E^2
\end{array} \right.
\end{aligned}
\end{align*} 
\end{tiny}
\end{theorem}
\subsubsection{"Box $\cap$ Polyhedral" Uncertainty}
When the uncertain set in the problem \eqref{eq4.7} is the intersection of a unit box and a unit polyhedron, it is represented by,
\begin{align}
\begin{aligned}	\label{eq4.30}
\mathscr{U}_{\mu}=\left\{\bm{\mu} : \bm{\mu}= \bm{\mu}^{(0)}+\sum_{j=1}^{n}\zeta_j \bm{\mu}^{(j)},\ ||\bm{\zeta}||_{\infty} \leq 1, \ ||\bm{\zeta}||_1 \leq 1 \right\}
\end{aligned}
\end{align}
Then any solution $\bm{x}$ is robust feasible for the uncertain Markowitz potfolio model with "box $\cap$ polyhedral" uncertainty, if and only if we have,
\begin{align}
\begin{aligned}	\label{eq4.31}
& \bm{x}^\top \left[\bm{\mu}^{(0)}+\sum_{j=1}^{n}\zeta_j \bm{\mu}^{(j)}\right] \left[\bm{\mu}^{(0)}+\sum_{j=1}^{n}\zeta_j \bm{\mu}^{(j)}\right]^\top\bm{x} - \tau^2 \geq 0, \qquad \forall \bm{\zeta}: ||\bm{\zeta}||_{\infty} \leq 1, \ ||\bm{\zeta}||_1 \leq 1
\end{aligned}
\end{align}
First we need to find the intersection between box and polyhedral uncertain sets of radii $\delta_B$ and $\delta_P$ respectively.\\
Now the box uncertain set of radius $\delta_B$ is given by,
\begin{align*}
\left\{\bm{\zeta} : ||\bm{\zeta}||_{\infty} \leq \delta_B \right\} \ \implies \left\{\zeta_j:|\zeta_j|\leq \delta_B, \quad \forall j\right\}
\end{align*}
Further, the polyhedral uncertain set of radius $\delta_P$ is given by,
\begin{align*}
\left\{\bm{\zeta} : ||\bm{\zeta}||_{1} \leq \delta_P\right\} \ \implies \left\{\zeta_j: \sum_{j=1}^{n} |\zeta_j| \leq \delta_P \right\}
\end{align*}
So we can write the combined uncertainty as intersection between the two sets. This can be written as,
\begin{align*}
\left\{\zeta_j: \sum_{j=1}^{n} |\zeta_j| \leq \min\left\{n\delta_B, \delta_P\right\}\right\} \ \implies \left\{\bm{\zeta}: ||\bm{\zeta}||_{1} \leq \min\left\{n\delta_B, \delta_P\right\} \right\} 
\end{align*}
Now there arises two cases.\\
\underline{\textit{\textbf{Case I:} When $\min\left\{n\delta_B, \delta_P\right\}=n\delta_B$}}\\
For this case the intersection set satisfies $||\bm{\zeta}||_{1} \leq n\delta_B$, and it can be achieved by making a transformation $\zeta_j \to \dfrac{\zeta_j}{n \delta_B}$ in the constraint \eqref{eq4.31}, which then becomes,
\begin{align}
\begin{aligned}	\label{eq4.32}
& \bm{x}^\top \left[\bm{\mu}^{(0)}(n\delta_B)+\sum_{j=1}^{n}\zeta_j \bm{\mu}^{(j)}\right] \left[\bm{\mu}^{(0)}(n\delta_B)+\sum_{j=1}^{n}\zeta_j \bm{\mu}^{(j)}\right]^\top\bm{x} - \tau^2 n^2\delta_B^2 \geq 0, \quad \forall \bm{\zeta}: ||\bm{\zeta}||_1 \leq n\delta_B
\end{aligned}
\end{align}
Equivalently the above condition is written as,
\begin{align}
\begin{aligned} \label{eq4.33}
\bm{x}^\top \left[\bm{\mu}^{(0)}(n\delta_B)+\sum_{j=1}^{n}\zeta_j \bm{\mu}^{(j)}\right] \left[\bm{\mu}^{(0)}(n\delta_B)+\sum_{j=1}^{n}\zeta_j \bm{\mu}^{(j)}\right]^\top\bm{x} - \tau^2 n^2\delta_B^2 \geq 0, \quad
\forall \bm{\zeta}: n^2\delta_B^2-\bm{\zeta}^\top S \bm{\zeta} \geq 0
\end{aligned}
\end{align}
where $S$ is an $(n \times n)$ matrix with $1$ as each entry given by \eqref{eq4.16c}.\\
The left hand side of the left expression of \eqref{eq4.33} into quadratic form $Q_{\bm x}^{(BP)}(\delta_{B}, \bm{\zeta})$, which is given by,
\begin{align}
\begin{aligned}	\label{eq4.33a}
&(\bm z_B)^\top \bm A'' (\bm z_B) \geq 0,
\end{aligned}
\end{align}
where  \begin{align}
\begin{aligned} \label{eq4.33b}
&\bm z_B=\begin{bmatrix} 
\delta_B\\
\zeta_1\\
\zeta_2\\
\vdots\\
\zeta_n
\end{bmatrix}\text{ is an $(n+1)$ dimensional vector, and }\\
&\bm A''=\begin{bmatrix} 
(n{\bm{\mu}^{(0)}}^\top\bm x)^2-n^2\tau^2 & (n{\bm{\mu}^{(0)}}^\top\bm x)({\bm{\mu}^{(1)}}^\top\bm x)  & (n{\bm{\mu}^{(0)}}^\top\bm x)({\bm{\mu}^{(2)}}^\top\bm x) & \dots & (n{\bm{\mu}^{(0)}}^\top\bm x)({\bm{\mu}^{(n)}}^\top\bm x)\\
({\bm{\mu}^{(1)}}^\top\bm x)(n{\bm{\mu}^{(0)}}^\top\bm x) & ({\bm{\mu}^{(1)}}^\top\bm x)^2 & ({\bm{\mu}^{(1)}}^\top\bm x)({\bm{\mu}^{(2)}}^\top\bm x) & \dots & ({\bm{\mu}^{(1)}}^\top\bm x)({\bm{\mu}^{(n)}}^\top\bm x)\\
({\bm{\mu}^{(2)}}^\top\bm x)(n{\bm{\mu}^{(0)}}^\top\bm x) & ({\bm{\mu}^{(2)}}^\top\bm x)({\bm{\mu}^{(1)}}^\top\bm x) & ({\bm{\mu}^{(2)}}^\top\bm x)^2 & \dots & ({\bm{\mu}^{(2)}}^\top\bm x)({\bm{\mu}^{(n)}}^\top\bm x)\\
\vdots & \vdots & \vdots & \ddots & \vdots \\
({\bm{\mu}^{(n)}}^\top\bm x)(n{\bm{\mu}^{(0)}}^\top\bm x) & ({\bm{\mu}^{(n)}}^\top\bm x)({\bm{\mu}^{(1)}}^\top\bm x) &({\bm{\mu}^{(n)}}^\top\bm x)({\bm{\mu}^{(2)}}^\top\bm x) & \dots & ({\bm{\mu}^{(n)}}^\top\bm x)^2
\end{bmatrix}
\end{aligned}
\end{align}
is a symmetric matrix of order $(n+1) \times (n+1)$.\\
Now writing the left hand side of the right expression of \eqref{eq4.33} in quadratic form, we have,
\begin{align}
\begin{aligned}	\label{eq4.33c}
& n^2\delta_B^2-\bm{\zeta}^\top S \bm{\zeta} &&= 
 {\begin{bmatrix} 
	\delta_B\\
	\zeta_1\\
	\zeta_2\\
	\vdots\\
	\zeta_n
	\end{bmatrix}}^\top
\begin{bmatrix} 
n^2 & 0  & 0 & \dots & 0\\
0 & -1 & -1 & \dots & -1\\
0 & -1 & -1 & \dots & -1\\
\vdots & \vdots & \vdots & \ddots & \vdots \\
0 & -1 & -1 & \dots & -1
\end{bmatrix}
\begin{bmatrix} 
\delta_B\\
\zeta_1\\
\zeta_2\\
\vdots\\
\zeta_n
\end{bmatrix}
\end{aligned}
\end{align}
Let us denote this quadratic form as, $P_1^{(BP)}(\delta_B, \bm{\zeta})$ and the right expression in \eqref{eq4.33} can be written as,
\begin{align}
\begin{aligned}	\label{eq4.33d}
& && {\begin{bmatrix} 
	\delta_B\\
	\zeta_1\\
	\zeta_2\\
	\vdots\\
	\zeta_n
	\end{bmatrix}}^\top
\begin{bmatrix} 
n^2 & 0  & 0 & \dots & 0\\
0 & -1 & -1 & \dots & -1\\
0 & -1 & -1 & \dots & -1\\
\vdots & \vdots & \vdots & \ddots & \vdots \\
0 & -1 & -1 & \dots & -1
\end{bmatrix}
\begin{bmatrix} 
\delta_B\\
\zeta_1\\
\zeta_2\\
\vdots\\
\zeta_n
\end{bmatrix} \geq 0
\end{aligned}
\end{align}
That means,
\begin{align}
\begin{aligned}	\label{eq4.33e}
&(\bm z_B)^\top \bm B_1^{(BP)} (\bm z_B) \geq 0,
\end{aligned}
\end{align}
where 
\begin{align} 
\begin{aligned}	\label{eq4.33f}
\bm B_1^{(BP)}=\begin{bmatrix} 
n^2 & 0  & 0 & \dots & 0\\
0 & -1 & -1 & \dots & -1\\
0 & -1 & -1 & \dots & -1\\
\vdots & \vdots & \vdots & \ddots & \vdots \\
0 & -1 & -1 & \dots & -1
\end{bmatrix}
\end{aligned}
\end{align}
is a symmetric matrix of order $(n+1) \times (n+1)$.\\
Now using the expressions \eqref{eq4.33a} and \eqref{eq4.33e}, we can write \eqref{eq4.33} as,
\begin{align}
\begin{aligned}	\label{eq4.34}
& (\bm z_B)^\top \bm A'' (\bm z_B) \geq 0, \quad \forall (\bm z_B) : (\bm z_B)^\top \bm B_1^{(BP)} (\bm z_B) \geq 0
\end{aligned}
\end{align}
Then by Lemma \ref{lemma4.1} we have, the matrix $\bm A'' - \lambda \bm B_1^{(BP)}$ is positive semidefinite.
Thus the following matrix inequality holds.
\begin{align}
\begin{aligned} \label{eq4.35}
\begin{scriptsize}
\begin{bmatrix} 
(n{\bm{\mu}^{(0)}}^\top\bm x)^2-n^2\tau^2-n^2\lambda & (n{\bm{\mu}^{(0)}}^\top\bm x)({\bm{\mu}^{(1)}}^\top\bm x)  & (n{\bm{\mu}^{(0)}}^\top\bm x)({\bm{\mu}^{(2)}}^\top\bm x) & \dots & (n{\bm{\mu}^{(0)}}^\top\bm x)({\bm{\mu}^{(n)}}^\top\bm x)\\
({\bm{\mu}^{(1)}}^\top\bm x)(n{\bm{\mu}^{(0)}}^\top\bm x) & ({\bm{\mu}^{(1)}}^\top\bm x)^2+\lambda & ({\bm{\mu}^{(1)}}^\top\bm x)({\bm{\mu}^{(2)}}^\top\bm x)+\lambda & \dots & ({\bm{\mu}^{(1)}}^\top\bm x)({\bm{\mu}^{(n)}}^\top\bm x)+\lambda\\
({\bm{\mu}^{(2)}}^\top\bm x)(n{\bm{\mu}^{(0)}}^\top\bm x) & ({\bm{\mu}^{(2)}}^\top\bm x)({\bm{\mu}^{(1)}}^\top\bm x)+\lambda & ({\bm{\mu}^{(2)}}^\top\bm x)^2+\lambda & \dots & ({\bm{\mu}^{(2)}}^\top\bm x)({\bm{\mu}^{(n)}}^\top\bm x)+\lambda\\
\vdots & \vdots & \vdots & \ddots & \vdots \\
({\bm{\mu}^{(n)}}^\top\bm x)(n{\bm{\mu}^{(0)}}^\top\bm x) & ({\bm{\mu}^{(n)}}^\top\bm x)({\bm{\mu}^{(1)}}^\top\bm x)+\lambda &({\bm{\mu}^{(n)}}^\top\bm x)({\bm{\mu}^{(2)}}^\top\bm x)+\lambda & \dots & ({\bm{\mu}^{(n)}}^\top\bm x)^2+\lambda
\end{bmatrix}
\end{scriptsize}
\succeq \bm 0
\end{aligned}
\end{align}
\underline{\textit{\textbf{Case II:} When $\min\left\{n\delta_B, \delta_P\right\}=\delta_P$}}\\
For this case the intersection set satisfies $||\bm{\zeta}||_{1} \leq \delta_P$, and it can be achieved by making a transformation $\zeta_j \to \dfrac{\zeta_j}{\delta_P}$ in the constraint \eqref{eq4.31}, which then becomes,
\begin{align}
\begin{aligned}	\label{eq4.36}
& \bm{x}^\top \left[\bm{\mu}^{(0)}(\delta_P)+\sum_{j=1}^{n}\zeta_j \bm{\mu}^{(j)}\right] \left[\bm{\mu}^{(0)}(\delta_P)+\sum_{j=1}^{n}\zeta_j \bm{\mu}^{(j)}\right]^\top\bm{x} - \tau^2 \delta_P^2 \geq 0, \qquad \forall \bm{\zeta}: ||\bm{\zeta}||_1 \leq \delta_P
\end{aligned}
\end{align}
This is equivalent to the polyhedral uncertainty case given in Theorem \ref{thm4.3}. So proceeding in a similar manner, the following matrix inequality holds.
\begin{align}
\begin{aligned} \label{eq4.36a}
\begin{scriptsize}
\begin{bmatrix} 
({\bm{\mu}^{(0)}}^\top\bm x)^2-\tau^2-\lambda & ({\bm{\mu}^{(0)}}^\top\bm x)({\bm{\mu}^{(1)}}^\top\bm x)  & ({\bm{\mu}^{(0)}}^\top\bm x)({\bm{\mu}^{(2)}}^\top\bm x) & \dots & ({\bm{\mu}^{(0)}}^\top\bm x)({\bm{\mu}^{(n)}}^\top\bm x)\\
({\bm{\mu}^{(1)}}^\top\bm x)({\bm{\mu}^{(0)}}^\top\bm x) & ({\bm{\mu}^{(1)}}^\top\bm x)^2+\lambda & ({\bm{\mu}^{(1)}}^\top\bm x)({\bm{\mu}^{(2)}}^\top\bm x)+\lambda & \dots & ({\bm{\mu}^{(1)}}^\top\bm x)({\bm{\mu}^{(n)}}^\top\bm x)+\lambda\\
({\bm{\mu}^{(2)}}^\top\bm x)({\bm{\mu}^{(0)}}^\top\bm x) & ({\bm{\mu}^{(2)}}^\top\bm x)({\bm{\mu}^{(1)}}^\top\bm x)+\lambda & ({\bm{\mu}^{(2)}}^\top\bm x)^2+\lambda & \dots & ({\bm{\mu}^{(2)}}^\top\bm x)({\bm{\mu}^{(n)}}^\top\bm x)+\lambda\\
\vdots & \vdots & \vdots & \ddots & \vdots \\
({\bm{\mu}^{(n)}}^\top\bm x)({\bm{\mu}^{(0)}}^\top\bm x) & ({\bm{\mu}^{(n)}}^\top\bm x)({\bm{\mu}^{(1)}}^\top\bm x)+\lambda &({\bm{\mu}^{(n)}}^\top\bm x)({\bm{\mu}^{(2)}}^\top\bm x)+\lambda & \dots & ({\bm{\mu}^{(n)}}^\top\bm x)^2+\lambda
\end{bmatrix}
\end{scriptsize}
\succeq \bm 0
\end{aligned}
\end{align} 
The results of the above two cases can be given as follows.
\begin{theorem} \label{thm4.5}
Consider the uncertain portfolio problem \eqref{eq4.5} under the uncertainty set obtained from the combination of box and polyhedral uncertainty sets. Let the radii of the box and polyhedral sets are respectively $\delta_B$ and $\delta_P$. Then the solution vector $\bm x$ to this problem is feasible for the robust constraint in \eqref{eq4.5} if and only if there exists a $\lambda \geq 0$ such that the following matrix inequalities hold:
\begin{tiny}
\begin{align*}
\begin{aligned}	
& \left \{ \begin{array}{ll}
\begin{bmatrix} 
(n{\bm{\mu}^{(0)}}^\top\bm x)^2-n^2\tau^2-n^2\lambda & (n{\bm{\mu}^{(0)}}^\top\bm x)({\bm{\mu}^{(1)}}^\top\bm x)  & (n{\bm{\mu}^{(0)}}^\top\bm x)({\bm{\mu}^{(2)}}^\top\bm x) & \dots & (n{\bm{\mu}^{(0)}}^\top\bm x)({\bm{\mu}^{(n)}}^\top\bm x)\\
({\bm{\mu}^{(1)}}^\top\bm x)(n{\bm{\mu}^{(0)}}^\top\bm x) & ({\bm{\mu}^{(1)}}^\top\bm x)^2+\lambda & ({\bm{\mu}^{(1)}}^\top\bm x)({\bm{\mu}^{(2)}}^\top\bm x)+\lambda & \dots & ({\bm{\mu}^{(1)}}^\top\bm x)({\bm{\mu}^{(n)}}^\top\bm x)+\lambda\\
({\bm{\mu}^{(2)}}^\top\bm x)(n{\bm{\mu}^{(0)}}^\top\bm x) & ({\bm{\mu}^{(2)}}^\top\bm x)({\bm{\mu}^{(1)}}^\top\bm x)+\lambda & ({\bm{\mu}^{(2)}}^\top\bm x)^2+\lambda & \dots & ({\bm{\mu}^{(2)}}^\top\bm x)({\bm{\mu}^{(n)}}^\top\bm x)+\lambda\\
\vdots & \vdots & \vdots & \ddots & \vdots \\
({\bm{\mu}^{(n)}}^\top\bm x)(n{\bm{\mu}^{(0)}}^\top\bm x) & ({\bm{\mu}^{(n)}}^\top\bm x)({\bm{\mu}^{(1)}}^\top\bm x)+\lambda &({\bm{\mu}^{(n)}}^\top\bm x)({\bm{\mu}^{(2)}}^\top\bm x)+\lambda & \dots & ({\bm{\mu}^{(n)}}^\top\bm x)^2+\lambda
\end{bmatrix} \succeq 0, &\mbox{ if } \min\left\{n\delta_B, \delta_P\right\}=n\delta_B\\
\\
\begin{bmatrix} 
({\bm{\mu}^{(0)}}^\top\bm x)^2-\tau^2-\lambda & ({\bm{\mu}^{(0)}}^\top\bm x)({\bm{\mu}^{(1)}}^\top\bm x)  & ({\bm{\mu}^{(0)}}^\top\bm x)({\bm{\mu}^{(2)}}^\top\bm x) & \dots & ({\bm{\mu}^{(0)}}^\top\bm x)({\bm{\mu}^{(n)}}^\top\bm x)\\
({\bm{\mu}^{(1)}}^\top\bm x)({\bm{\mu}^{(0)}}^\top\bm x) & ({\bm{\mu}^{(1)}}^\top\bm x)^2+\lambda & ({\bm{\mu}^{(1)}}^\top\bm x)({\bm{\mu}^{(2)}}^\top\bm x)+\lambda & \dots & ({\bm{\mu}^{(1)}}^\top\bm x)({\bm{\mu}^{(n)}}^\top\bm x)+\lambda\\
({\bm{\mu}^{(2)}}^\top\bm x)({\bm{\mu}^{(0)}}^\top\bm x) & ({\bm{\mu}^{(2)}}^\top\bm x)({\bm{\mu}^{(1)}}^\top\bm x)+\lambda & ({\bm{\mu}^{(2)}}^\top\bm x)^2+\lambda & \dots & ({\bm{\mu}^{(2)}}^\top\bm x)({\bm{\mu}^{(n)}}^\top\bm x)+\lambda\\
\vdots & \vdots & \vdots & \ddots & \vdots \\
({\bm{\mu}^{(n)}}^\top\bm x)({\bm{\mu}^{(0)}}^\top\bm x) & ({\bm{\mu}^{(n)}}^\top\bm x)({\bm{\mu}^{(1)}}^\top\bm x)+\lambda &({\bm{\mu}^{(n)}}^\top\bm x)({\bm{\mu}^{(2)}}^\top\bm x)+\lambda & \dots & ({\bm{\mu}^{(n)}}^\top\bm x)^2+\lambda\\
\end{bmatrix} \succeq 0, &\mbox{ if } \min\left\{n\delta_B, \delta_P\right\}=\delta_P
\end{array} \right.
\end{aligned}
\end{align*} 
\end{tiny}
\end{theorem}

\subsubsection{"Ellipsoidal $\cap$ Polyhedral" Uncertainty}
When the uncertain set in the problem \eqref{eq4.7} is the intersection of a unit box and a unit polyhedron, it is represented by,
\begin{align}
\begin{aligned}	\label{eq4.43}
\mathscr{U}_{\mu}=\left\{\bm{\mu} : \bm{\mu}= \bm{\mu}^{(0)}+\sum_{j=1}^{n}\zeta_j \bm{\mu}^{(j)},\ ||\bm{\zeta}||_1 \leq 1, \ ||\bm{\zeta}||_2 \leq 1 \right\}
\end{aligned}
\end{align}
%where $||\cdot||_1$ represents a polyhedral set and $||\cdot||_2$ represents an ellipsoidal set.\\
Then any solution $\bm{x}$ is robust feasible for the uncertain Markowitz potfolio model with "ellipsoidal $\cap$ polyhedral" uncertainty, if and only if we have,
\begin{align}
\begin{aligned}	\label{eq4.44}
& \bm{x}^\top- \left[\bm{\mu}^{(0)}+\sum_{j=1}^{n}\zeta_j \bm{\mu}^{(j)}\right] \left[\bm{\mu}^{(0)}+\sum_{j=1}^{n}\zeta_j \bm{\mu}^{(j)}\right]^\top\bm{x} - \tau^2 \geq 0, \qquad \forall \bm{\zeta}: ||\bm{\zeta}||_1 \leq 1, \ ||\bm{\zeta}||_2 \leq 1
\end{aligned}
\end{align}
First we need to find the intersection between ellipsoidal and polyhedral uncertain sets of radii $\delta_E$ and $\delta_P$ respectively.\\
Now the ellipsoidal uncertain set of radius $\delta_E$ is given by,
\begin{align*}
	\left\{\bm{\zeta} : ||\bm{\zeta}||_{2} \leq \delta_E \right\} \implies \left\{\zeta_j:\sum_{j=1}^{n}(\zeta_j)^2 \leq \delta_E^2 \right\}
\end{align*}
Further, the polyhedral uncertain set of radius $\delta_P$ is given by,
\begin{align*}
	\left\{\bm{\zeta} : ||\bm{\zeta}||_{1} \leq \delta_P\right\}  \implies \left\{\zeta_j: \left(\sum_{j=1}^{n} |\zeta_j| \right)^2 \leq \delta_P^2 \right\}
\end{align*}
So we can write the combined uncertainty as intersection between the two sets. This can be written as,
\begin{align*}
\begin{aligned}
& &&\left\{\zeta_j: \sum_{j=1}^{n} (\zeta_j)^2 \leq \min\{\delta_E^2, \delta_P^2-2\sum_{i \neq j}{|\zeta_i| |\zeta_j|}\}\right\} \\
&\implies && \left\{\zeta_j: \sum_{j=1}^{n} (\zeta_j)^2 \leq \min\{\delta_E^2, \delta_P^2\}\right\},
\end{aligned}
\end{align*} 
since each $|\zeta_j|\geq 0$, which implies $\delta_P^2-2\sum_{i \neq j}{|\zeta_i| |\zeta_j|} \leq \delta_P^2$.\\
So the combined uncertain set is given by,
\begin{align*}
\begin{aligned}
& && \left\{\bm{\zeta}: ||\bm{\zeta}||_{2} \leq \min\{\delta_E, \delta_P \} \right\}
\end{aligned}
\end{align*} 
Now there arises two cases.\\
\underline{\textit{\textbf{Case I:} When $\min\left\{\delta_E, \delta_P\right\}=\delta_E$}}\\
For this case the intersection set satisfies $||\bm{\zeta}||_{2} \leq \delta_E$, and it can be achieved by making a transformation $\zeta_j \to \dfrac{\zeta_j}{\delta_E}$ in the constraint \eqref{eq4.44}, which then becomes,
\begin{align}
\begin{aligned}	\label{eq4.45}
& \bm{x}^\top \left[\bm{\mu}^{(0)}(\delta_E)+\sum_{j=1}^{n}\zeta_j \bm{\mu}^{(j)}\right] \left[\bm{\mu}^{(0)}(\delta_E)+\sum_{j=1}^{n}\zeta_j \bm{\mu}^{(j)}\right]^\top\bm{x} - \tau^2\delta_E^2 \geq 0, \quad \forall \bm{\zeta}: ||\bm{\zeta}||_2 \leq \delta_E
\end{aligned}
\end{align}
Equivalently the above constraint is written as,
\begin{align}
\begin{aligned} \label{eq4.46}
& \bm{x}^\top \left[\bm{\mu}^{(0)}(\delta_E)+\sum_{j=1}^{n}\zeta_j \bm{\mu}^{(j)}\right] \left[\bm{\mu}^{(0)}(\delta_E)+\sum_{j=1}^{n}\zeta_j \bm{\mu}^{(j)}\right]^\top\bm{x} - \tau^2\delta_E^2 \geq 0, \  \forall \bm{\zeta}: \delta_E^2-\bm{\zeta}^\top \bm{\zeta} \geq 0
\end{aligned}
\end{align}
This is equivalent to the ellipsoidal uncertainty case given in Theorem \ref{thm4.1}. So proceeding in a similar manner, the following matrix inequality holds.
\begin{align}
\begin{aligned} \label{eq4.46a}
\begin{scriptsize}
\begin{bmatrix} 
({\bm{\mu}^{(0)}}^\top\bm x)^2-\tau^2-\lambda & ({\bm{\mu}^{(0)}}^\top\bm x)({\bm{\mu}^{(1)}}^\top\bm x)  & ({\bm{\mu}^{(0)}}^\top\bm x)({\bm{\mu}^{(2)}}^\top\bm x) & \dots & ({\bm{\mu}^{(0)}}^\top\bm x)({\bm{\mu}^{(n)}}^\top\bm x)\\
({\bm{\mu}^{(1)}}^\top\bm x)({\bm{\mu}^{(0)}}^\top\bm x) & ({\bm{\mu}^{(1)}}^\top\bm x)^2+\lambda & ({\bm{\mu}^{(1)}}^\top\bm x)({\bm{\mu}^{(2)}}^\top\bm x) & \dots & ({\bm{\mu}^{(1)}}^\top\bm x)({\bm{\mu}^{(n)}}^\top\bm x)\\
({\bm{\mu}^{(2)}}^\top\bm x)({\bm{\mu}^{(0)}}^\top\bm x) & ({\bm{\mu}^{(2)}}^\top\bm x)({\bm{\mu}^{(1)}}^\top\bm x)& ({\bm{\mu}^{(2)}}^\top\bm x)^2+\lambda & \dots & ({\bm{\mu}^{(2)}}^\top\bm x)({\bm{\mu}^{(n)}}^\top\bm x)\\
\vdots & \vdots & \vdots & \ddots & \vdots \\
({\bm{\mu}^{(n)}}^\top\bm x)({\bm{\mu}^{(0)}}^\top\bm x) & ({\bm{\mu}^{(n)}}^\top\bm x)({\bm{\mu}^{(1)}}^\top\bm x) &({\bm{\mu}^{(n)}}^\top\bm x)({\bm{\mu}^{(2)}}^\top\bm x) & \dots & ({\bm{\mu}^{(n)}}^\top\bm x)^2+\lambda
\end{bmatrix}
\end{scriptsize}
\succeq \bm 0
\end{aligned}
\end{align}
\underline{\textit{\textbf{Case II:} When $\min\left\{\delta_E, \delta_P\right\}=\delta_P$}}\\
For this case the intersection set satisfies $||\bm{\zeta}||_{2} \leq \delta_P$, and it can be achieved by making a transformation $\zeta_j \to \dfrac{\zeta_j}{\delta_P}$ in the constraint \eqref{eq4.44}, which then becomes,
\begin{align}
\begin{aligned}	\label{eq4.47}
& \bm{x}^\top \left[\bm{\mu}^{(0)}(\delta_P)+\sum_{j=1}^{n}\zeta_j \bm{\mu}^{(j)}\right] \left[\bm{\mu}^{(0)}(\delta_P)+\sum_{j=1}^{n}\zeta_j \bm{\mu}^{(j)}\right]^\top\bm{x} - \tau^2 \delta_P^2 \geq 0, \quad \forall \bm{\zeta}: ||\bm{\zeta}||_2 \leq \delta_P\\
\end{aligned}
\end{align}
Equivalently the above constraint is written as,
\begin{align}
\begin{aligned} \label{eq4.48}
& \bm{x}^\top \left[\bm{\mu}^{(0)}(\delta_P)+\sum_{j=1}^{n}\zeta_j \bm{\mu}^{(j)}\right] \left[\bm{\mu}^{(0)}(\delta_P)+\sum_{j=1}^{n}\zeta_j \bm{\mu}^{(j)}\right]^\top\bm{x} - \tau^2\delta_P^2 \geq 0, \  \forall \bm{\zeta}: \delta_P^2-\bm{\zeta}^\top \bm{\zeta} \geq 0
\end{aligned}
\end{align}
This is equivalent to the case of ellipsoidal uncertain set with radius $\delta_P$. So the similar derivations as Theorem \ref{thm4.1} reduces the constraint \eqref{eq4.48} to,
\begin{align}
\begin{aligned}	\label{eq4.49}
& (\bm z_P)^\top \bm A \bm (z_P) \geq 0, \quad \forall (\bm z_P) : (\bm z_P)^\top \bm B_2^{(EP)} (\bm z_P) \geq 0
\end{aligned}
\end{align}
where $z_P$, $A$, and $B_2^{(EP)}$ are respectively given by,
\begin{align}
\begin{aligned} \label{eq4.50}
& \bm z_P=\begin{bmatrix} 
\delta_P\\
\zeta_1\\
\zeta_2\\
\vdots\\
\zeta_n
\end{bmatrix}, \ \bm A=\begin{bmatrix} 
({\bm{\mu}^{(0)}}^\top\bm x)^2-\tau^2 & ({\bm{\mu}^{(0)}}^\top\bm x)({\bm{\mu}^{(1)}}^\top\bm x)  & ({\bm{\mu}^{(0)}}^\top\bm x)({\bm{\mu}^{(2)}}^\top\bm x) & \dots & ({\bm{\mu}^{(0)}}^\top\bm x)({\bm{\mu}^{(n)}}^\top\bm x)\\
({\bm{\mu}^{(1)}}^\top\bm x)({\bm{\mu}^{(0)}}^\top\bm x) & ({\bm{\mu}^{(1)}}^\top\bm x)^2 & ({\bm{\mu}^{(1)}}^\top\bm x)({\bm{\mu}^{(2)}}^\top\bm x) & \dots & ({\bm{\mu}^{(1)}}^\top\bm x)({\bm{\mu}^{(n)}}^\top\bm x)\\
({\bm{\mu}^{(2)}}^\top\bm x)({\bm{\mu}^{(0)}}^\top\bm x) & ({\bm{\mu}^{(2)}}^\top\bm x)({\bm{\mu}^{(1)}}^\top\bm x) & ({\bm{\mu}^{(2)}}^\top\bm x)^2 & \dots & ({\bm{\mu}^{(2)}}^\top\bm x)({\bm{\mu}^{(n)}}^\top\bm x)\\
\vdots & \vdots & \vdots & \ddots & \vdots \\
({\bm{\mu}^{(n)}}^\top\bm x)({\bm{\mu}^{(0)}}^\top\bm x) & ({\bm{\mu}^{(n)}}^\top\bm x)({\bm{\mu}^{(1)}}^\top\bm x) &({\bm{\mu}^{(n)}}^\top\bm x)({\bm{\mu}^{(2)}}^\top\bm x) & \dots & ({\bm{\mu}^{(n)}}^\top\bm x)^2
\end{bmatrix}\\
& \bm B^{(EP)}=\begin{bmatrix} 
1 & 0  & 0 & \dots & 0\\
0 & -1 & 0 & \dots & 0\\
0 & 0 & -1 & \dots & 0\\
\vdots & \vdots & \vdots & \ddots & \vdots \\
0 & 0 &0 & \dots & -1
\end{bmatrix}
\end{aligned}
\end{align}
Then using Lemma \ref{lemma4.1} we have, the matrix $\bm A - \lambda \bm B^{(EP)}$ is positive semidefinite. Thus the following matrix inequality holds.

\begin{align}
\begin{aligned} \label{eq4.51}
\begin{scriptsize}
\begin{bmatrix} 
({\bm{\mu}^{(0)}}^\top\bm x)^2-\tau^2-\lambda & ({\bm{\mu}^{(0)}}^\top\bm x)({\bm{\mu}^{(1)}}^\top\bm x)  & ({\bm{\mu}^{(0)}}^\top\bm x)({\bm{\mu}^{(2)}}^\top\bm x) & \dots & ({\bm{\mu}^{(0)}}^\top\bm x)({\bm{\mu}^{(n)}}^\top\bm x)\\
({\bm{\mu}^{(1)}}^\top\bm x)({\bm{\mu}^{(0)}}^\top\bm x) & ({\bm{\mu}^{(1)}}^\top\bm x)^2+\lambda & ({\bm{\mu}^{(1)}}^\top\bm x)({\bm{\mu}^{(2)}}^\top\bm x) & \dots & ({\bm{\mu}^{(1)}}^\top\bm x)({\bm{\mu}^{(n)}}^\top\bm x)\\
({\bm{\mu}^{(2)}}^\top\bm x)({\bm{\mu}^{(0)}}^\top\bm x) & ({\bm{\mu}^{(2)}}^\top\bm x)({\bm{\mu}^{(1)}}^\top\bm x)& ({\bm{\mu}^{(2)}}^\top\bm x)^2+\lambda & \dots & ({\bm{\mu}^{(2)}}^\top\bm x)({\bm{\mu}^{(n)}}^\top\bm x)\\
\vdots & \vdots & \vdots & \ddots & \vdots \\
({\bm{\mu}^{(n)}}^\top\bm x)({\bm{\mu}^{(0)}}^\top\bm x) & ({\bm{\mu}^{(n)}}^\top\bm x)({\bm{\mu}^{(1)}}^\top\bm x) &({\bm{\mu}^{(n)}}^\top\bm x)({\bm{\mu}^{(2)}}^\top\bm x) & \dots & ({\bm{\mu}^{(n)}}^\top\bm x)^2+\lambda
\end{bmatrix}
\end{scriptsize}
\succeq \bm 0
\end{aligned}
\end{align}
Therefore the feasibility condition over "ellipsoidal $\cap$ polyhedral" uncertainty remains the same for both the cases and is given by the following result.
\begin{theorem} \label{thm4.6}
Consider the uncertain portfolio problem \eqref{eq4.5} under the uncertainty set obtained from the combination of ellipsoidal and polyhedral uncertainty sets. Let the radii of the ellipsoidal and polyhedral sets are respectively $\delta_E$ and $\delta_P$. Then the solution vector $\bm x$ to this problem is feasible for the robust constraint in \eqref{eq4.5} if and only if there exists a $\lambda \geq 0$ such that the following matrix inequality hold:
\begin{align*}
\begin{scriptsize}
\begin{bmatrix} 
({\bm{\mu}^{(0)}}^\top\bm x)^2-\tau^2-\lambda & ({\bm{\mu}^{(0)}}^\top\bm x)({\bm{\mu}^{(1)}}^\top\bm x)  & ({\bm{\mu}^{(0)}}^\top\bm x)({\bm{\mu}^{(2)}}^\top\bm x) & \dots & ({\bm{\mu}^{(0)}}^\top\bm x)({\bm{\mu}^{(n)}}^\top\bm x)\\
({\bm{\mu}^{(1)}}^\top\bm x)({\bm{\mu}^{(0)}}^\top\bm x) & ({\bm{\mu}^{(1)}}^\top\bm x)^2+\lambda & ({\bm{\mu}^{(1)}}^\top\bm x)({\bm{\mu}^{(2)}}^\top\bm x) & \dots & ({\bm{\mu}^{(1)}}^\top\bm x)({\bm{\mu}^{(n)}}^\top\bm x)\\
({\bm{\mu}^{(2)}}^\top\bm x)({\bm{\mu}^{(0)}}^\top\bm x) & ({\bm{\mu}^{(2)}}^\top\bm x)({\bm{\mu}^{(1)}}^\top\bm x)& ({\bm{\mu}^{(2)}}^\top\bm x)^2+\lambda & \dots & ({\bm{\mu}^{(2)}}^\top\bm x)({\bm{\mu}^{(n)}}^\top\bm x)\\
\vdots & \vdots & \vdots & \ddots & \vdots \\
({\bm{\mu}^{(n)}}^\top\bm x)({\bm{\mu}^{(0)}}^\top\bm x) & ({\bm{\mu}^{(n)}}^\top\bm x)({\bm{\mu}^{(1)}}^\top\bm x) &({\bm{\mu}^{(n)}}^\top\bm x)({\bm{\mu}^{(2)}}^\top\bm x) & \dots & ({\bm{\mu}^{(n)}}^\top\bm x)^2+\lambda
\end{bmatrix}
\end{scriptsize} \succeq 0.
\end{align*}
\end{theorem}

\section{Conclusion} \label{sec4.4}
\par In this paper, we have obtained the feasibility conditions for the robust counterparts of the uncertain Markowitz model in the form of nonlinear semidefinite programming. We represented our uncertainty sets in the form of quadratic functions and that helped us to write our uncertain constraint in quadratic forms of the perturbation vector $\bm{\zeta}$. We can also express the robust counterpart problem in form of a quadratic semidefinite programming problem by replacing the uncertain constraint with the obtained feasibility condition. In this way, we will also be able to obtain the robust counterparts for the combined uncertainty sets, which in general is difficult to find out. With this approach, we can make use of the solution algorithms, convergence analysis for non-linear semidefinite programming given in \cite{Correa 2004,Freund 2007}.
\par In addition, the study of combined uncertainty sets is more important because the robust solutions under these sets are less conservative as compared to the individual uncertainty sets. There is a scope of further research in obtaining the robust counterparts of some other portfolio models in the semidefinite programming form by using this approach.
\section*{Declaration on Conflict of Interest} The authors hereby declare that they have no conflicts of interest.
%\begin{acknowledgements}
%If you'd like to thank anyone, place your comments here
%and remove the percent signs.
%\end{acknowledgements}

% Authors must disclose all relationships or interests that 
% could have direct or potential influence or impart bias on 
% the work: 
%
% \section*{Conflict of interest}
%
% The authors declare that they have no conflict of interest.

% BibTeX users please use one of
%\bibliographystyle{spbasic}      % basic style, author-year citations
%\bibliographystyle{spmpsci}      % mathematics and physical sciences
%\bibliographystyle{spphys}       % APS-like style for physics
%\bibliography{}   % name your BibTeX data base

% Non-BibTeX users please use

\end{document}